\documentclass[10pt]{article}
\usepackage[utf8]{inputenc}
\usepackage[T1]{fontenc}
\usepackage{amsmath}
\usepackage{amsfonts}
\usepackage{amssymb}
\usepackage[version=4]{mhchem}
\usepackage{stmaryrd}
\usepackage{hyperref}
\hypersetup{colorlinks=true, linkcolor=blue, filecolor=magenta, urlcolor=cyan,}
\usepackage{bbold}
\usepackage{graphicx}
\usepackage[export]{adjustbox}
\graphicspath{ {./images/} }

\title{Stationary states of aggregation-diffusion equations with compactly supported attraction kernels: radial symmetry and mass-independent boundedness }

\author{Roumen Anguelov, Chelsea Bright\\
Department of Mathematics and Applied Mathematics\\
 University of Pretoria,  South Africa\\
roumen.anguelov@up.ac.za,\ chels.alex.bright@gmail.com}
\date{}

\let\svthefootnote\thefootnote
\newcommand\blfootnotetext[1]{%
  \let\thefootnote\relax\footnote{#1}%
  \addtocounter{footnote}{-1}%
  \let\thefootnote\svthefootnote%
}

\let\svfootnotetext\footnotetext
\renewcommand\footnotetext[2][?]{%
  \if\relax#1\relax%
    \ifnum\value{footnote}=0\blfootnotetext{#2}\else\svfootnotetext{#2}\fi%
  \else%
    \if?#1\ifnum\value{footnote}=0\blfootnotetext{#2}\else\svfootnotetext{#2}\fi%
    \else\svfootnotetext[#1]{#2}\fi%
  \fi
}

\begin{document}
\maketitle

\begin{abstract}
We consider a nonlocal aggregation diffusion equation incorporating repulsion modelled by nonlinear diffusion and attraction modelled by nonlocal interaction. When the attractive interaction kernel is radially symmetric and strictly increasing on its domain it is previously known that all stationary solutions are radially symmetric and decreasing up to a translation; however, this result has not been extended to accommodate attractive kernels that are non-decreasing, for instance, attractive kernels with bounded support. For the diffusion coefficient $m>1$, we show that, for attractive kernels that are radially symmetric and non-decreasing, all stationary states are radially symmetric and decreasing up to a translation on each connected subset of their support. Furthermore, for $m>2$, we prove analytically that stationary states have an upper-bound independent of the initial data, confirming previous numerical results given in the literature.
\end{abstract}


\section*{1. Introduction}
General aggregation diffusion equations have been used in a variety of different applications, including the modelling of chemotaxis, the biological aggregation of insects and herds of animals, see $[6,26,27,30,1]$, and opinion dynamics, see [12] and the references therein. The behaviour in these settings are typically driven by long-range attraction and short-range repulsion. The competition between these forces gives rise to characteristic time-independent morphologies. These resulting equilibria are the focus of this study. In particular, we consider a continuum description of the collective behaviour where we
analyze the evolution of the population density $\rho(x, t)$ at some location $x \in \mathbb{R}^{d}$ and at time $t \geq 0$.

We consider the initial value problem given by the following non-local integro-differential equation

\begin{equation*}
\partial_{t} \rho=\varepsilon \Delta \rho^{m}+\nabla \cdot(\rho \nabla(W * \rho)) \quad x \in \mathbb{R}^{d}, t>0 \tag{1}
\end{equation*}

with initial condition $\rho_{0} \in L_{+}^{1}\left(\mathbb{R}^{d}\right) \cap L^{m}\left(\mathbb{R}^{d}\right)$. The local repulsion is modelled using nonlinear diffusion with $m>1$ and where $\varepsilon$ is the diffusion coefficient. The non-local attraction arises from the second term on the right, where $W * \rho=\int_{\mathbb{R}^{d}} W(x-y) \rho(y) d y$. In essence, the presence of individuals at position $y \in \mathbb{R}^{d}$ creates a force, proportional to $-\nabla W(x-y)$, that acts on the individuals positioned at $x \in \mathbb{R}^{d}$. The interaction kernel $W: \mathbb{R}^{d} \rightarrow \mathbb{R}$ is given and is assumed to be radially symmetric and non-decreasing from its centre. That is, there is a function $\omega:[0, \infty) \rightarrow \mathbb{R}$ such that $W(x)=\omega(\|x\|)$ for all $x \in \mathbb{R}^{d}$ and $\omega^{\prime}(x) \geq 0, x \in \mathbb{R}_{+}$. Thus, the interaction kernel only takes into account attractive effects.

Significant work has been done on the aggregation diffusion equation where the interaction kernel is the attractive power-law kernel of the form

$$
W_{k}(x)=\left\{\begin{array}{l}
\frac{|x|^{k}}{k}, k \neq 0 \\
\ln |x|, k=0
\end{array}\right.
$$

for $2-d \leq k \leq 2[17,5,4,2,8,15]$. A special case of a power-law kernel where $k=2-d$ is the Newtonian kernel. In fact, if the convolution operator $W *$ acting on $\rho$ in Equation (1) is the Newtonian kernel then we obtain the nonlinear parabolic elliptic Keller-Segel model of chemotaxis, see [7, 10]. This follows from the fact that the Newtonian potential is the inverse of the negative Laplacian.

Mathematically, the first question regarding equation (1) is the existence and uniqueness of solutions. This is almost completely resolved in [19]. More precisely, it is proved in [19] that, for the case of general attractive kernels that are no more singular at the origin than the Newtonian kernel and with either bounded or unbounded support, weak solutions of Problem (1) exist globally in time for any $m>2-\frac{2}{d}$.

In the sequel, we assume that these conditions hold. The assumptions on $W$ made so far can be given as follows:\\
$W 1 W(x) \in C^{1}\left(\mathbb{R}^{d} \backslash\{0\}\right)$ and is radially symmetric. That is, there exists a function $\omega:(0, \infty) \rightarrow \mathbb{R}$ such that $W(x)=\omega(\|x\|)$.\\
$W 2$ There exists some $C_{\omega}>0$ such that $\omega^{\prime}(r) \leq C_{\omega} r^{1-d}$ for $r \leq 1$.

The solutions that exist satisfy the property that $\int_{\mathbb{R}^{d}} \rho(x, t) d x=\int_{\mathbb{R}^{d}} \rho_{0}(x) d x$, that is, the total mass of the solution is conserved for all time. Additionally, $\rho(x, t) \geq 0$ for all $t \geq 0[11]$.

Considering weak solutions of Problem (1), we have that $\rho_{s} \in L_{+}^{1}\left(\mathbb{R}^{d}\right)$ is a stationary state of Equation (1) if $\rho_{s}^{m} \in H_{l o c}^{1}\left(\mathbb{R}^{d}\right), \nabla W * \rho_{s} \in L_{l o c}^{1}\left(\mathbb{R}^{d}\right)$, and it satisfies

\begin{equation*}
\varepsilon \nabla \rho_{s}^{m}=-\rho_{s}\left(\nabla W * \rho_{s}\right) \text { in } \operatorname{supp} \rho_{s} \tag{2}
\end{equation*}

in the sense of distributions in $\mathbb{R}^{d}$.\\
Many of the results in the literature regarding existence and stability of stationary states rely on the fact that there is a Lyapunov functional for the dynamics of Equation (1) given by the energy functional,

\begin{equation*}
\mathcal{E}[\rho]=\frac{\varepsilon}{m-1} \int_{\mathbb{R}^{d}} \rho^{m} d x+\frac{1}{2} \int_{\mathbb{R}^{d}} \rho(W * \rho) d x=: \mathcal{S}[\rho]+\mathcal{I}[\rho] \tag{3}
\end{equation*}

A general method to prove existence of stationary states of (1) is by showing that the global minimizer of (3) corresponds to a stationary state of (1). It is proved in [3] and [17] that existence of a global minimizer depends on the choice of $m$ and $W$. In particular, if $W$ is an integrable and attractive kernel, the author of [3] proves that $m=2$ is the threshold value of $m$ such that for $m>2$, existence holds for any given mass, while for $1<m \leq 2$, the existence depends on the value of $\|W\|_{L^{1}\left(\mathbb{R}^{d}\right)}$. Under more general assumptions on the interaction kernel $W$, it is proved in [17] that in the diffusion dominated regime, that is $m>\max \left\{2-\frac{2}{d}, 1\right\}$, a global minimizer exists for any given mass.

By Reisz's rearrangement inequality it is known that a global minimizer of (3) must be radially decreasing [17]. Hence, it is natural to ask whether the same is true for stationary states of (1). It is shown in [17] that for any $m>0$, all stationary solutions of (1) are indeed radially symmetric and decreasing up to a translation for the case where $\omega^{\prime}$ has unbounded support. It is further proved in [20] that, for $m \geq 2$, the stationary state is unique up to translation. However, to the best of our knowledge, proof of the radial symmetry and uniqueness of stationary states for the case where the support of $\omega^{\prime}$ is bounded has not been given.

Our first main result in this work is that, for $m>1$ and the case where $\omega^{\prime}$ has bounded support, a stationary solution of (1) is radially symmetric and decreasing up to a translation on each connected component of its support.

Our second main result is that for $m>2$ the stationary states have an upper-bound independent of the initial data for both cases of $\omega^{\prime}$ with bounded and unbounded support. This is a notable characteristic of the model as it is a natural property of physical and biological aggregations, where individuals will aggregate together up to a maximum density and no further. The value of the\\
diffusion coefficient $m$ plays a key role in the emergence of this boundedness, where $m=2$ is a threshold. More precisely, we prove analytically that for $m>2$ we obtain boundedness of stationary states independent of the initial data, while for $m \leq 2$ the maximal density grows with the total mass $\int_{\mathbb{R}^{d}} \rho_{0}(x) d x$. This boundedness property for the special case of $m=3$ is shown numerically in [30].

Throughout this paper, we assume that, in addition to ( $W 1$ ) and (W2), one of the following assumptions on the attractive kernel $W$ holds:\\
$W 3 \omega^{\prime}(r)>0$ for all $r>0$ and there exists some $C_{\omega}>0$ such that $\omega^{\prime}(r) \leq C_{\omega}$ for $r>1$. Moreover, $\lim _{r \rightarrow \infty} \omega(r)=0$ and there exists an $\alpha \in(0, d)$ for which $m>1+\frac{\alpha}{d}$ and $\omega(\tau r) \leq \tau^{-\alpha} \omega(r)$ for all $\tau \geq 1$ and $r>0$.\\
$\hat{W} 4$ There exists $q>0$ such that $\omega^{\prime}(r)>0$ for all $0<r<q$ and $\omega^{\prime}(r)=0$ for all $r \geq q$.

Assumption (W3) combines assumptions $(K 3),(K 4)$, and (K6) given in [17], allowing us to refer to the results obtained therein. Note that in [17] it is assumed that $\lim _{r \rightarrow \infty} \omega(r)=\ell \in(0, \infty)$. However, it is not restrictive to assume that $\ell=0$ since adding a constant to the potential $W$ does not change Equation $(1)$.

Note that the threshold value $r=1$ in (W2) and $(W 3)$ can be replaced by any positive value since, in essence, $(W 2)$ places restriction on $\omega^{\prime}$ when $r$ is small, while ( $W 3$ ) places a restriction on $\omega^{\prime}$ when $r$ is large. In addition, the threshold value $r=q$ in $(\hat{W} 4)$ can be assumed to be 1. Therefore, for simplicity, we reformulate ( $\hat{W} 4)$ as follows:\\
$W 4 \omega^{\prime}(r)>0$ for all $0<r<1$ and $\omega^{\prime}(r)=\omega(r)=0$ for all $r \geq 1$.\\
Our main consideration in this paper is when $W$ satisfies assumption (W4). That is, when $\omega^{\prime}$ has compact support. With reference to our first main result, we recall that the radial symmetry property for $\omega^{\prime}$ with unbounded support is proved in [17]. In terms of the notations adopted here, the results in [17, Theorem 2.2, Theorem 3.1, Theorem 3.7, Lemma 3.8, Lemma 3.9] can be formulated as follows:

Theorem 1.1. [17] For $m>\max \left\{2-\frac{2}{d}, 1\right\}$, assume that conditions (W1), (W2), and (W3) hold. Then, for any positive mass $M$, there exists a stationary state $\rho_{s} \in L_{+}^{1}\left(\mathbb{R}^{d}\right) \cap L^{\infty}\left(\mathbb{R}^{d}\right)$ of (1). Furthermore, any stationary state in $L_{+}^{1}\left(\mathbb{R}^{d}\right) \cap L^{\infty}\left(\mathbb{R}^{d}\right)$ is radially symmetric and decreasing up to $a$ translation and compactly supported.

In Section 2 we extend this result in an appropriate way for attractive kernels with compact support, that is, when (W4) holds. More precisely, we show that stationary states are radially symmetric and decreasing up to a translation on each connected component of their support. Moreover, we prove that if $\operatorname{supp} \rho_{s}$ has more than one connected component, then the distance between any two\\
components is not smaller than the radius of supp $W$.\\
Further, the structure of the paper is as follows. The mass-independent boundedness of stationary states for $m>2$ is proved in Section 3. In Section 4, we provide numerical results visualizing the radial symmetry and boundedness properties of the stationary states of (1). Some concluding remarks are given in Section 5 .

\section*{2. Stationary states for compactly supported attractive kernels}
In this section, we consider stationary solutions of Equation (1) where the attractive kernel $W$ has compact support. We prove, for this case of attractive kernel, that continuous, compactly supported stationary solutions of (1) exist. Our main result of this section is that, for $m>1$, stationary states are radially decreasing up to a translation on each connected component of their support, where the definition of a connected component is given by Definition 2.3. Furthermore, we prove that, for $m>2$, if the support of a stationary state has more than one connected component then the distance between any two components is at least the radius of supp $W$.

We use the following definitions in the results of this section.\\
Definition 2.1. Radial symmetry of a non-negative function\\
A non-negative function $f$ on $\mathbb{R}^{d}$ is radially symmetric if there is a function $\tilde{f}$ defined on $[0, \infty)$ such that $f(x)=\tilde{f}(\|x\|)$ for all $x \in \mathbb{R}^{d}$.

Definition 2.2. Radially decreasing up to a translation\\
A non-negative function $f$ on $\mathbb{R}^{d}$ is radially decreasing up to a translation if there exists some $x_{0} \in \mathbb{R}^{d}$ such that $f\left(\cdot-x_{0}\right)$ is radially symmetric and $\tilde{f}\left(\left\|x-x_{0}\right\|\right)$ is non-increasing in $\left\|x-x_{0}\right\|$.

Here, $\| \cdot||$ denotes the Euclidean norm.\\
Definition 2.3 (Connected component). The connected components of a topological space $X$ are closed, disjoint, non-empty subsets of $X$ such that their union is the whole space $X$.

We summarize our main results of this section in the following theorem.\\
Theorem 2.1. Let $m>2$ and let $W$ satisfy assumptions (W1), (W2), and (W4). Then there exists a stationary solution $\rho_{s} \in L_{+}^{1}\left(\mathbb{R}^{d}\right) \cap L^{\infty}\left(\mathbb{R}^{d}\right)$ of (1). Furthermore, $\rho_{s} \in C\left(\mathbb{R}^{d}\right)$, and is radially symmetric, decreasing, and compactly supported on each connected component of $\operatorname{supp} \rho_{s}$. Additionally, if $\operatorname{supp} \rho_{s}$ has more than one connected component then the distance between any two components is at least the radius of $\operatorname{supp} W$.

\subsection*{2.1. Existence of stationary states}
As previously mentioned, global minimizers of the energy functional can be used to determine existence of stationary states of (1). Hence, we start by stating a previous result on the existence of radially decreasing global minimizers of (3).

Theorem 2.2. [3] For $m>2, W \in L^{1}\left(\mathbb{R}^{d}\right)$ radially symmetric and nondecreasing, and for any $M>0$, there exists a radially symmetric and decreasing global minimizer of the energy functional (3) defined in

$$
\mathcal{Y}_{M}:=\left\{\rho \in L_{+}^{1}\left(\mathbb{R}^{d}\right) \cap L^{m}\left(\mathbb{R}^{d}\right):\|\rho\|_{L^{1}\left(\mathbb{R}^{d}\right)}=M\right\}
$$

It is easy to see that $W \in L^{1}\left(\mathbb{R}^{d}\right)$ under assumptions $(W 1)$, ( $W 2$ ), and $(W 4)$. Hence, we have existence of a radially symmetric and decreasing global minimizer of (3). We further show, under these assumptions on the attractive kernel, that all global minimizers of (3) defined in $\mathcal{Y}_{M}$, whose support consists of a single connected component, are radially symmetric and decreasing, compactly supported, uniformly bounded, and correspond to stationary states of $(1)$ in the weak sense.

Theorem 2.3. Let $m>2$ and let $W$ satisfy assumptions (W1), (W2), and (W4). If $\bar{\rho}$ is a global minimizer of (3) in $\mathcal{Y}_{M}$ whose support consists of a single connected component, then $\bar{\rho}$ is radially symmetric and decreasing up to a translation.

Proof. The proof follows similarly to the proof of [13, Proposition 4.3], where it is assumed that $\bar{\rho}$ is not radially symmetric and decreasing under any translation. Using the fact that $-W(x)=-\omega(\|x\|)$ is non-negative, radially symmetric, and non-increasing on $\|x\|>0$, Riesz's rearrangement inequality is then used to obtain a contradiction.

The proof of Theorem 2.4 below follows similarly to the proof of [18, Theorem 3.1], which considers global minimizers of the energy functional corresponding to the $2 D$ Keller-Segel equation, that is, where $W$ is the Newtonian kernel. Theorem 2.4 is used to show that global minimizers of the energy functional are in fact stationary states of (1).

Theorem 2.4. Let $m>2$ and let $W$ satisfy assumptions ( $W 1$ ), (W2), and (W4). If $\bar{\rho}$ is a global minimizer of (3) in $\mathcal{Y}_{M}$, then there exists a constant $D[\bar{\rho}]$ such that

\begin{equation*}
\frac{m \varepsilon}{m-1} \bar{\rho}^{m-1}(x)+(W * \bar{\rho})(x)=D[\bar{\rho}] \text {, a.e in } \operatorname{supp} \bar{\rho} \tag{4}
\end{equation*}

and

\begin{equation*}
\frac{m \varepsilon}{m-1} \bar{\rho}^{m-1}(x)+(W * \bar{\rho})(x) \geq D[\bar{\rho}] \text {, a.e outside } \operatorname{supp} \bar{\rho} \tag{5}
\end{equation*}

where

$$
D[\bar{\rho}]=\frac{2}{M} \mathcal{E}[\bar{\rho}]+\frac{m-2}{M(m-1)}\|\bar{\rho}\|_{m}^{m}
$$

Using Theorems 2.3 and 2.4, we are able to obtain the following result:\\
Theorem 2.5. Let $m>2$ and let $W$ satisfy assumptions (W1), (W2), and (W4). If $\bar{\rho}$ is a global minimizer of (3) in $\mathcal{Y}_{M}$ whose support consists of a single connected component, then $\bar{\rho}$ is compactly supported.

Proof. Suppose for a contradiction that supp $\bar{\rho}$ is not compact. Then supp $\bar{\rho}=$ $\mathbb{R}^{d}$ since, by Theorem $2.3, \bar{\rho}$ is radially decreasing from some $x_{0} \in \mathbb{R}^{d}$ acting as a centre. Hence, from (4), we have that there exists a constant $C$ such that

\begin{equation*}
\frac{m \varepsilon}{m-1} \bar{\rho}^{m-1}(x)+(W * \bar{\rho})(x)=C \tag{6}
\end{equation*}

for a.e. $x \in \mathbb{R}^{d}$. Since $\bar{\rho}$ is radially decreasing and in $L^{1}\left(\mathbb{R}^{d}\right)$, there is a function $\bar{\rho}_{*}$ where $\bar{\rho}(x)=\bar{\rho}_{*}(|| x||)$ and where $\lim _{\|x\| \rightarrow \infty} \bar{\rho}_{*}(\|x\|)=0$. In addition, we claim that $\lim _{\|x\| \rightarrow \infty}(W * \bar{\rho})(x)=0$. To show this, let $\mathcal{A}=\left\{y \in \mathbb{R}^{d}: \| x-\right.$ $y \|<1\}$ and fix $\|x\|>\left\|x_{0}\right\|+1$. Then, since $\|x-y\|<1$ implies $\left\|x_{0}\right\|<$ $\|x\|-1<\|y\|$, we have that $\bar{\rho}_{*}\left(\left\|x_{0}\right\|\right) \geq \bar{\rho}_{*}(\|x\|-1) \geq \bar{\rho}_{*}(\|y\|)$. Therefore, for $\|x\|>\left\|x_{0}\right\|+1$

$$
\int_{\mathbb{R}^{d}} W(x-y) \bar{\rho}(y) d y=\int_{A} W(x-y) \bar{\rho}_{*}(\|y\|) d y \leq \bar{\rho}_{*}(\|x\|-1)\|W\|_{L^{1}\left(\mathbb{R}^{d}\right)}
$$

Taking the limit as $\|x\| \rightarrow \infty$ yields our claim. Now, taking the limit as $\|x\| \rightarrow \infty$ in (6), we obtain

$$
\frac{m \varepsilon}{m-1} \bar{\rho}^{m-1}(x)=(-W * \bar{\rho})(x)
$$

for a.e. $x \in \mathbb{R}^{d}$. Fix $x \in \mathbb{R}^{d}$ and define $\mathcal{B}=\{y \in \mathcal{A}:\|x\|>\|y\|\}$. Since $W(x)=\omega(\|x\|) \leq 0$ for all $\|x\|>0$, we have that $|W(x)|=-W(x)$ for all $x \in \mathbb{R}^{d}$ and so there exists $C>0$ such that

\begin{equation*}
\int_{\mathbb{R}^{d}}-W(x-y) \bar{\rho}(y) d y \geq \bar{\rho}_{*}(\|x\|) \int_{\mathcal{B}}|W(x-y)| d y=C \bar{\rho}_{*}(\| x||)>0 \tag{7}
\end{equation*}

Thus, $\frac{m \varepsilon}{m-1} \bar{\rho}^{m-1}(x) \geq C \bar{\rho}(x)>0$, yielding $\bar{\rho}(x) \geq\left(\frac{C(m-1)}{m \varepsilon}\right)^{\frac{1}{m-2}}>0$, for a.e $x \in \mathbb{R}^{d}$. Hence,

$$
\int_{\mathbb{R}^{d}} \bar{\rho}(x) d x \geq \int_{\mathbb{R}^{d}}\left(\frac{C(m-1)}{m \varepsilon}\right)^{\frac{1}{m-2}} d x=\infty
$$

This contradicts the fact that $\bar{\rho} \in L^{1}\left(\mathbb{R}^{d}\right)$. Hence, supp $\bar{\rho}$ is compact.

Theorem 2.6. Let $m>2$ and let $W$ satisfy assumptions (W1), (W2), and (W4) and let $\bar{\rho} \in \mathcal{Y}_{M}$ be a global minimizer of (3) whose support consists of a single connected component. Then $\bar{\rho} \in L^{\infty}\left(\mathbb{R}^{d}\right)$.

Theorems 2.3, 2.4, and 2.5 are applied to obtain the above result. The method of proof for Theorem 2.6 follows similarly to that of [17, Lemma 3.9], where it is proved that the boundedness property of $\bar{\rho}$ holds under the assumptions ( $W 1$ )-(W3). Indeed, if ( $W 4$ ) holds, then $\lim _{r \rightarrow \infty} \omega(r)=0$ and there exists $C_{\omega}>0$ such that $\omega^{\prime}(r) \leq C_{\omega}$ for all $r>1$, as given in ( $\left.W 3\right)$. Furthermore, the proof of [17, Lemma 3.9] does not rely on the assumptions made in (W3) that $\omega^{\prime}(r)>0$ for all $r>0$ and that there exists an $\alpha \in(0, d)$ for which $m>1+\frac{\alpha}{d}$ and $\omega(\tau r) \leq \tau^{-\alpha} \omega(r)$ for all $\tau \geq 1$ and $r>0$.

Using Theorem 2.6 above, it is easy to see that $W * \bar{\rho} \in \mathcal{W}^{1, \infty}\left(\mathbb{R}^{d}\right)$. Indeed, fixing $x \in \mathbb{R}^{d}$ and letting $\mathcal{A}=\left\{y \in \mathbb{R}^{d}:\|x-y\|<1\right\}$, we have that

\begin{align*}
|(\nabla W * \bar{\rho})(x)| & \leq \int_{\mathbb{R}^{d}} \omega^{\prime}(\|x-y\|) \bar{\rho}(y) d y \\
& \leq C_{\omega}\|\bar{\rho}\|_{L^{\infty}\left(\mathbb{R}^{d}\right)} \int_{\mathcal{A}} \frac{1}{\|x-y\|^{d-1}} d y=C:=\text { Const. } \tag{8}
\end{align*}

Furthermore,\\
$|(W * \bar{\rho})(x)| \leq \int_{\mathbb{R}^{d}}|W(x-y)| \bar{\rho}(y) d y \leq\|\bar{\rho}\|_{L^{\infty}\left(\mathbb{R}^{d}\right)}\|W\|_{L^{1}\left(\mathbb{R}^{d}\right)}$, for all $x \in \mathbb{R}^{d}$\\
We now have the results necessary to relate global minimizers of the energy functional to stationary solutions of Equation (1).

Theorem 2.7. Let $m>2$ and let $W$ satisfy assumptions (W1), (W2), and (W4) and let $\bar{\rho} \in \mathcal{Y}_{M}$ be a global minimizer of (3) whose support consists of a single connected component. Then $\bar{\rho}$ is a stationary solution of (1) in the weak sense.

Proof. From Theorem 2.4 we have that

\begin{equation*}
\frac{m \varepsilon}{m-1} \bar{\rho}^{m-1}+(W * \bar{\rho})=D[\bar{\rho}], \text { a.e in } \operatorname{supp} \bar{\rho} \tag{10}
\end{equation*}

Furthermore, since $W * \bar{\rho} \in \mathcal{W}^{1, \infty}\left(\mathbb{R}^{d}\right)$, we can take gradients on both sides of (10) and multiply by $\bar{\rho}$ to obtain

$$
\frac{m \varepsilon}{m-1} \bar{\rho} \nabla \bar{\rho}^{m-1}=-\bar{\rho} \nabla(W * \bar{\rho}), \text { a.e in } \operatorname{supp} \bar{\rho}
$$

Using the fact that $\bar{\rho} \nabla \bar{\rho}^{m-1}=\frac{m-1}{m} \nabla \bar{\rho}^{m}$, we have that

$$
\varepsilon \nabla \bar{\rho}^{m}=-\bar{\rho} \nabla(W * \bar{\rho}), \text { a.e in } \operatorname{supp} \bar{\rho}
$$

yielding (2), as required.\\
From Theorem 2.7 we have the existence of a stationary solution $\rho_{s}$ of (1) in $L_{+}^{1}\left(\mathbb{R}^{d}\right) \cap L^{\infty}\left(\mathbb{R}^{d}\right)$. We now show that, for any stationary state of (1) in\\
$L_{+}^{1}\left(\mathbb{R}^{d}\right) \cap L^{\infty}\left(\mathbb{R}^{d}\right)$ with support possibly made up of more than one connected component, we have that on each connected component of $\operatorname{supp} \rho_{s}, \rho_{s}$ is radially symmetric and decreasing up to a translation and compactly supported. Additionally, the distance between any two components is at least the radius of $\operatorname{supp} W$.

In order to prove our result on radial symmetry of stationary solutions for compactly supported $W$, given in the next section, we need Lemma 2.8 given below. A similar result, for the case when $(W 1)-(W 3)$ are satisfied, is given in [17, Lemma 2.3]. Note that in [17, Lemma 2.3] an extra assumption is made, that is, $\omega(1+\|x\|) \rho_{s} \in L^{1}\left(\mathbb{R}^{d}\right)$. This assumption need not be made for the case of $\omega$ satisfying (W1),(W2), and (W4). This is because $\omega(r)=0$ for all $r \geq 1$, by assumption $(W 4)$, which yields $\omega(1+\|x\|) \rho_{s} \in L^{1}\left(\mathbb{R}^{d}\right)$, as required.

Lemma 2.8. Let $\rho_{s} \in L_{+}^{1}\left(\mathbb{R}^{d}\right) \cap L^{\infty}\left(\mathbb{R}^{d}\right)$ be a non-negative stationary state of (1) where $m>1$ and where $W$ satisfies assumptions ( $W 1$ ), (W2), and (W4). Then $\rho_{s} \in C\left(\mathbb{R}^{d}\right)$ and

\begin{equation*}
\frac{m \varepsilon}{m-1} \rho_{s}^{m-1}(x)+\left(W * \rho_{s}\right)(x)=C_{j}, \text { for } x \in D_{j} \tag{11}
\end{equation*}

where $C_{j}$ may be different on each connected component $D_{j}$ of $\operatorname{supp} \rho_{s}$.\\
Furthermore, there exists some $C=C\left(\left\|\rho_{s}\right\|_{L^{1}},\left\|\rho_{s}\right\|_{L^{\infty}}, C_{\omega}, d\right)>0$ such that

$$
\frac{m \varepsilon}{m-1}\left|\nabla\left(\rho_{s}^{m-1}\right)\right| \leq C \text { in } \operatorname{supp} \rho_{s}
$$

Proof. The proof follows the same approach as that of [17, Lemma 2.3]. It is only necessary to prove, under the assumptions $(W 1),(W 2)$, and (W4), that $\nabla W * \rho_{s}$ and $W * \rho_{s}$ are globally bounded. This is proved by using the same argument given to obtain (8) and (9).

From Theorem 2.7, we have that global minimizers of the energy, whose supports are connected, are stationary states of Equation (1). The converse is not necessarily true; however, it is the case that a stationary state of (1) is a stationary point of the energy, as stated in the following theorem.

Theorem 2.9. If $\rho_{s}$ is a stationary state of Equation (1), then $\rho_{s}$ is a stationary point of the energy functional (3).

Theorem 2.9 is proved using the same approach as the proof of Equation (3.25) in [18, Theorem 3.1]. That is, for the same choice of test functions, it can be shown that the first variation of the energy exists and vanishes under the assumption that $\rho_{s}$ is a stationary state of (1) and, thus, satisfies (11).

\subsection*{2.2. Radial symmetry of stationary states}
In this subsection, we prove our first main result, namely that, for $m>1$ and $W$ satisfying assumptions $(W 1),(W 2)$, and ( $W 4$ ), all stationary solutions of (1) are radially symmetric and decreasing when restricted to a single connected component of their support. Furthermore, we show that for $m>2$ and any stationary state whose support has more than one connected component, the distance between any two components is at least the radius of the support of $W$. These results are summarized in the theorem below.

Theorem 2.10. Let $W$ satisfy assumptions ( $W 1$ ), (W2), and (W4). Let $\rho_{s} \in L_{+}^{1}\left(\mathbb{R}^{d}\right) \cap L^{\infty}\left(\mathbb{R}^{d}\right)$ be a stationary state of (1) and let $D \subset \mathbb{R}^{d}$ be a connected component of $\operatorname{supp} \rho_{s}$, as given in Definition 2.3. That is, $D \cap$ closure $\left(\operatorname{supp} \rho_{s} \backslash D\right)=\emptyset$. Then the following holds:

\begin{enumerate}
  \item For any $m>1$, there exists an $x_{0} \in D$ such that $\left.\rho_{s}\right|_{D}$ is radially symmetric and decreasing from $x_{0}$ as a centre.
  \item For any $m>2$, for all $x \in \operatorname{interior}(D)$ and for any $y \in \operatorname{interior}\left(\operatorname{supp} \rho_{s} \backslash\right.$ $D)$, we have $\|x-y\| \geq 1$.
\end{enumerate}

Proof outline for Theorem 2.10.1\\
We prove Theorem 2.10 .1 by contradiction, assuming $\rho_{s}$ is not radially symmetric and decreasing under any translation. Similar to the proof of [17, Theorem 2.2], we use continuous Steiner symmetrization to construct a family of densities $\mu(\tau, \cdot)$ with $\left.\rho_{s}\right|_{D}=\mu(0, \cdot)=: \mu_{0}$ such that $\mathcal{E}[\mu(\tau)]-\mathcal{E}\left[\mu_{0}\right]<-C_{1} \tau$, for some $C_{1}>0$ and any sufficiently small $\tau>0$. However; since $\mu_{0}$ is a stationary state of $\mathcal{E}$, it can be shown that $\left|\mathcal{E}[\mu(\tau)]-\mathcal{E}\left[\mu_{0}\right]\right| \leq C_{2} \tau^{2}$ for some $C_{2}>0$ and for all $\tau$ small enough. Combining these two inequalities results in a contradiction for $\tau$ sufficiently small.

We first give the definition of continuous Steiner symmetrization of a set or function. We refer to [17], [22], and [25] for a deeper introduction.

Definition 2.4. The continuous Steiner symmetrization of any open set $U \subset \mathbb{R}$, denoted $M^{\tau}(U), \tau \geq 0$, is defined as below. We denote $I(\tilde{x}+c, r):=(\tilde{x}+c-$ $r, \tilde{x}+c+r)$, where $\tilde{x}, c, r \in \mathbb{R}$.

\begin{enumerate}
  \item If $U=I(\tilde{x}+c, r)$, then
\end{enumerate}

$$
M^{\tau}(I(\tilde{x}+c, r)):= \begin{cases}I(\tilde{x}+c-\tau \operatorname{sgn} c, r), & \text { if } 0 \leq \tau<|c| \\ I(\tilde{x}, r), & \text { if } \tau \geq|c|\end{cases}
$$

\begin{enumerate}
  \setcounter{enumi}{1}
  \item If $U=\cup_{i=1}^{N} I\left(\tilde{x}+c_{i}, r_{i}\right)$, where all $I\left(\tilde{x}+c_{i}, r_{i}\right)$ are disjoint from each other, then $M^{\tau}(U):=\cup_{i=1}^{N} M^{\tau}\left(I\left(\tilde{x}+c_{i}, r_{i}\right)\right)$ for $0 \leq \tau<\tau_{1}$, where $\tau_{1}$ is the first time two intervals $M^{\tau}\left(I\left(\tilde{x}+c_{i}, r_{i}\right)\right)$ share a common endpoint. Once this occurs, we merge the two intervals sharing a common endpoint into one open interval and then define $M^{\tau}(U)$ in the same way, starting from $\tau=\tau_{1}$.\\
    \item If $U=\cup_{i=1}^{\infty} I\left(\tilde{x}+c_{i}, r_{i}\right)$, let $U_{N}=\cup_{i=1}^{N} I\left(\tilde{x}+c_{i}, r_{i}\right)$, for each $N \geq 1$, and define $M^{\tau}(U):=\cup_{N=1}^{\infty} M^{\tau}\left(U_{N}\right)$.
\end{enumerate}

We define the continuous Steiner symmetrization of a non-negative function $f$ defined on $\mathbb{R}$ as follows:

Definition 2.5. Let $f \in L_{+}^{1}(\mathbb{R})$. The continuous Steiner symmetrization for $f$ is given by

$$
S^{\tau} f(\alpha):=\int_{0}^{\infty} \chi_{M^{\tau}\left(U^{h}\right)}(\alpha) d h
$$

Similarly, we may define the continuous Steiner symmetrization of a nonnegative function $f$ defined on $\mathbb{R}^{d}$, with respect to a specific direction, as follows:\\
Definition 2.6. Fix $x_{1} \in \mathbb{R}, x^{\prime} \in \mathbb{R}^{d-1}$, and $h>0$. Let $f \in L_{+}^{1}\left(\mathbb{R}^{d}\right)$. The continuous Steiner symmetrization for $f$ in the direction $x_{1}$ is given by

$$
S^{\tau} f\left(x_{1}, x^{\prime}\right):=\int_{0}^{\infty} \chi_{M^{\tau}\left(U_{x^{\prime}}^{h}\right)}\left(x_{1}\right) d h
$$

It is easy to see that $\lim _{\tau \rightarrow \infty} S^{\tau} f=S f$, where $S f$ denotes the translation of the Steiner symmetrization $f$ such that $S f$ is symmetric decreasing about the hyperplane $\left\{x_{1}=\tilde{x}\right\}$.

In order to prove Theorem 2.10.1, we must prove the preliminary results, Lemma 2.11 and 2.12, given below, where we investigate how the interaction energy, denoted by $\mathcal{I}$ in (3), between two densities $\mu_{1}, \mu_{2} \in C(\mathbb{R})$ changes under their continuous Steiner symmetrizations. That is, we consider how

\begin{equation*}
I_{K}\left[\mu_{1}, \mu_{2}\right](\tau):=\int_{\mathbb{R} \times \mathbb{R}} S^{\tau} \mu_{1}(\alpha) S^{\tau} \mu_{2}(\beta) K(\alpha-\beta) d \alpha d \beta \tag{12}
\end{equation*}

changes with respect to $\tau$, given some $K \in C^{1}(\mathbb{R})$ to be defined later.\\
Lemma 2.11. Assume $K \in C^{1}(\mathbb{R})$ is an even function with $K^{\prime}(z)<0$ for all $0<z<R$ and $K^{\prime}(z)=0$ for all $|z| \geq R$. Let $\mu_{i}:=\chi_{I\left(\tilde{x}+c_{i}, r_{i}\right)}$ for $i=1,2$. Then, for $I(\tau):=I_{K}\left[\mu_{1}, \mu_{2}\right](\tau)$,

\begin{enumerate}
  \item $\frac{d^{+}}{d \tau} I(0) \geq 0$.
  \item In addition, if $\operatorname{sgn} c_{1} \neq \operatorname{sgn} c_{2}$,
\end{enumerate}

\begin{equation*}
\left|c_{2}-c_{1}\right|<r_{2}+r_{1}+R, \quad \text { and } \quad\left|r_{2}-r_{1}\right|<\left|c_{2}-c_{1}\right|+R \tag{13}
\end{equation*}

then

$$
\frac{d^{+}}{d \tau} I(0) \geq \frac{1}{6} \varphi\left(c_{1}, r_{1}, c_{2}, r_{2}, R\right) \min _{r \in\left[\frac{R}{3 \sqrt{2}}, \frac{R}{\sqrt{2}}\right]} K^{\prime}(r)=: c>0
$$

where

\begin{equation*}
\varphi\left(c_{1}, r_{1}, c_{2}, r_{2}, R\right)=\min \left\{-\left|r_{1}-r_{2}\right|+\left|c_{2}-c_{1}\right|+R,-\left|c_{2}-c_{1}\right|+r_{1}+r_{2}+R, R\right\} \tag{14}
\end{equation*}

Proof. Without loss of generality, assume $c_{2} \geq c_{1}$. For the case where $c_{2}<c_{1}$ the roles of $x$ and $y$ in the proof are reversed. For $x \in \mathbb{R}$, we have that $S^{\tau} \mu_{i}=$ $\int_{0}^{\infty} \chi_{M^{\tau}\left(U^{h}\left(\mu_{i}\right)\right)}(x) d h=\chi_{M^{\tau}\left(I\left(\tilde{x}+c_{i}, r_{i}\right)\right)}$ for $i=1,2$. Then for any $\tau \geq 0$ it follows that

$$
\begin{aligned}
I(\tau) & =I_{K}\left[\chi_{I\left(\tilde{x}+c_{1}, r_{1}\right)}, \chi_{I\left(\tilde{x}+c_{2}, r_{2}\right)}\right](\tau) \\
& =\int_{-r_{1}+\tilde{x}+c_{1}-\tau \operatorname{sgn} c_{1}}^{r_{1}+\tilde{x}+c_{1}-\tau \operatorname{sgn} c_{1}} \int_{-r_{2}+\tilde{x}+c_{2}-\tau \operatorname{sgn} c_{2}}^{r_{2}+\tilde{x}+c_{2}-\tau \operatorname{sgn} c_{2}} K(x-y) d y d x \\
& =\int_{-r_{1}}^{r_{1}} \int_{-r_{2}}^{r_{2}} K\left(x-y+\left(c_{1}-c_{2}\right)+\tau\left(\operatorname{sgn} c_{2}-\operatorname{sgn} c_{1}\right)\right) d y d x
\end{aligned}
$$

If $\operatorname{sgn} c_{1}=\operatorname{sgn} c_{2}$, then

$$
\frac{d^{+}}{d \tau} I(0)=\left(\operatorname{sgn} c_{2}-\operatorname{sgn} c_{1}\right) \int_{-r_{1}}^{r_{1}} \int_{-r_{2}}^{r_{2}} K^{\prime}\left(x-y+\left(c_{1}-c_{2}\right)\right) d y d x=0
$$

If $\operatorname{sgn} c_{1} \neq \operatorname{sgn} c_{2}$ we have that $\operatorname{sgn} c_{2}-\operatorname{sgn} c_{1}$ is either 2 or 1 . Hence,

\begin{equation*}
\frac{d^{+}}{d \tau} I(0)=\left(\operatorname{sgn} c_{2}-\operatorname{sgn} c_{1}\right) \int_{Q} K^{\prime}(x-y) d y d x \tag{15}
\end{equation*}

where $Q$ is the rectangle $\left[c_{1}-r_{1}, c_{1}+r_{1}\right] \times\left[c_{2}-r_{2}, c_{2}+r_{2}\right]$. Define $Q^{-}:=Q \cap\{0<$ $x-y<R\}$ and $Q^{+}:=Q \cap\{-R<x-y<0\}$. Note that $K^{\prime}(x-y)<0$ in $Q^{-}$ and $K^{\prime}(x-y)>0$ in $Q^{+}$. Since $K^{\prime}(x-y)=0$ for $|x-y| \geq R$, we have that

$$
\frac{d^{+}}{d \tau} I(0) \geq \int_{Q^{+}} K^{\prime}(x-y) d y d x+\int_{Q^{-}} K^{\prime}(x-y) d y d x
$$

Regardless of the choice of $r_{1}$ and $r_{2}$, since $\operatorname{sgn} c_{1} \neq \operatorname{sgn} c_{2}$ and $c_{2}>c_{1}, Q$ forms a rectangle with its centre, given by $\left(c_{1}, c_{2}\right)$, lying above the line $y=x$. Hence, for any $h>0$, the length of the line segment $Q^{+} \cap\{x-y=-h\}$ will be greater or equal to the length of $Q^{-} \cap\{x-y=h\}$. This implies that $\left|Q^{+}\right| \geq\left|Q^{-}\right|$and so $\frac{d^{+}}{d \tau} I(0) \geq 0$, which proves 1 .

Now, assume that $\operatorname{sgn} c_{1} \neq \operatorname{sgn} c_{2}$ and (13) holds. Furthermore, assume that $r_{2} \geq r_{1}$. Under these assumptions we obtain three possibilities, given in Figure 1 and 2 (a). We see in Figure 1 (a) that, since $r_{2}-r_{1}<\left|c_{2}-c_{1}\right|+R$, the bottom left-hand corner of the rectangle must be above the line $y=x-R$. Similarly, in Figure $1(\mathrm{~b})$, since $\left|c_{2}-c_{1}\right|<r_{2}+r_{1}+R$, we see that the bottom right-hand corner of the rectangle must be below the line $y=x+R$. Under our assumptions we will always have that the area of $Q^{+}$will be strictly greater than that\vspace{3mm}

\includegraphics[max width=\textwidth, center]{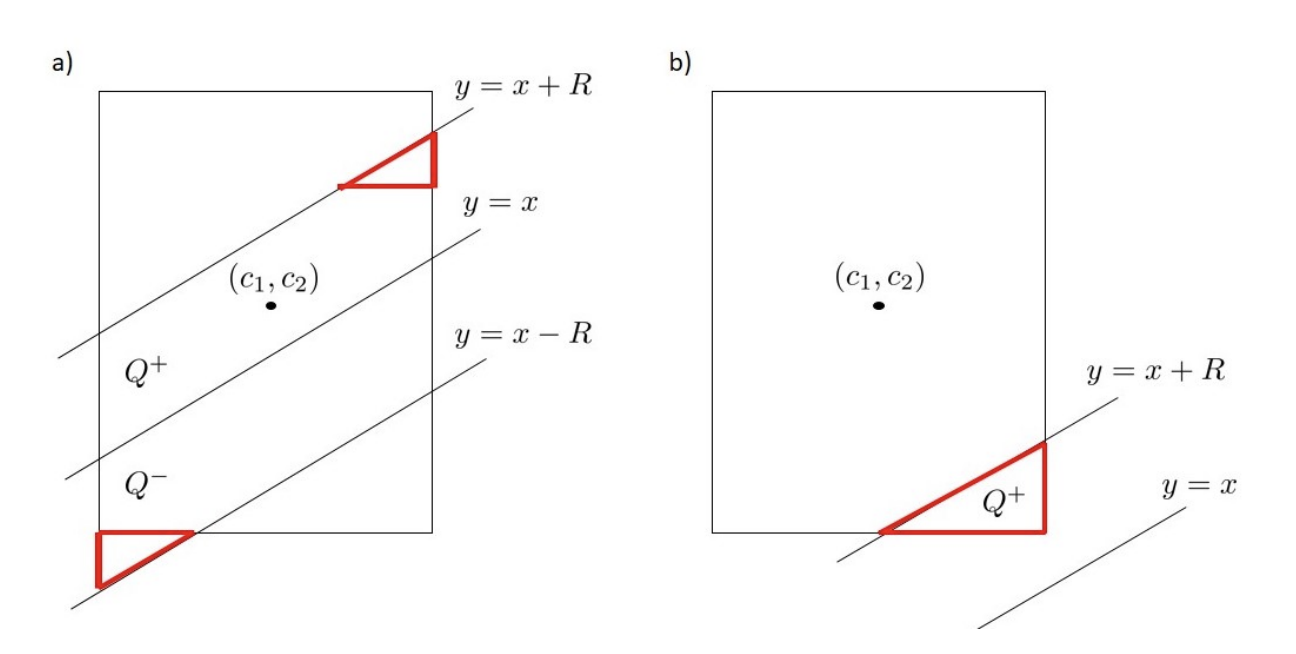}
Figure 1: Illustration of the rectangle $Q$ given in the proof of Lemma 2.11.\\
of $Q^{-}$, regardless of the choice of $r_{1}$ and $r_{2}$, where $r_{2} \geq r_{1}$. Furthermore, in all cases the difference in area will be at least the size of the triangle denoted $D$ in Figure 2 (b) (outlined red in Figure 1 (a), (b), and Figure 2 (a)).\vspace{3mm}

The vertices of $D$ are given by $\left(c_{1}+r_{1}, c_{1}+r_{1}+R\right),\left(c_{1}+r_{1}, z\right)$ and $(z-R, z)$, where $z:=\max \left\{2 c_{1}-c_{2}+r_{2}, c_{1}+r_{1}, c_{2}-r_{2}\right\}$. Now, consider the trapezium $\Omega \subset D$ where the bases of $\Omega$ lie parallel to the hypotenuse of $D$ and the longer base intersects the medicentre of $D$, as illustrated in Figure 2 (b). Since $K^{\prime}(x-$ $y)>0$ for all $x, y \in \Omega$, we have that

$$
\frac{d^{+}}{d \tau} I(0) \geq \int_{\Omega} K^{\prime}(x-y) d y d x \geq|\Omega| \min _{(x, y) \in \Omega} K^{\prime}(x-y)
$$

Now, substituting $\left(c_{1}+r_{1}, z\right)$ and the medicentre of $D$, given by $\frac{1}{3}\left(2 c_{1}+\right.$ $\left.2 r_{1}+z-R, c_{1}+r_{1}+R+2 z\right)$, into the equation of the unit normal to $y=x$, given by $\frac{1}{\sqrt{2}}(y-x)=0$, we find that $\min _{(x, y) \in \Omega} K^{\prime}(x-y)=\min _{r \in\left[z_{1}, z_{2}\right]} K^{\prime}(r)$, where

$$
z_{1}=\frac{1}{\sqrt{2}}\left(z-c_{1}-r_{1}\right) \geq \frac{R}{3 \sqrt{2}}
$$

and

$$
z_{2}=\frac{1}{3 \sqrt{2}}\left(z-c_{1}-r_{1}+2 R\right) \leq \frac{1}{3 \sqrt{2}}\left(c_{1}+r_{1}+R-c_{1}-r_{1}+2 R\right)=\frac{R}{\sqrt{2}}
$$

Hence,

$$
\min _{r \in\left[z_{1}, z_{2}\right]} K^{\prime}(r) \geq \min _{r \in\left[\frac{R}{3 \sqrt{2}}, \frac{R}{\sqrt{2}}\right]} K^{\prime}(r)>0
$$

Furthermore, by the properties of the medicentre, the longer base of $\Omega$ divides both the base and the height of $D$ in a ratio $1: 2$ as illustrated in Figure 2 (b). Therefore, by the definition of $\Omega$, we find that $|\Omega|=\frac{1}{6}\left(c_{1}+r_{1}-z+R\right)^{2}$.\\[3mm]
\begin{center}
\includegraphics[max width=\textwidth]{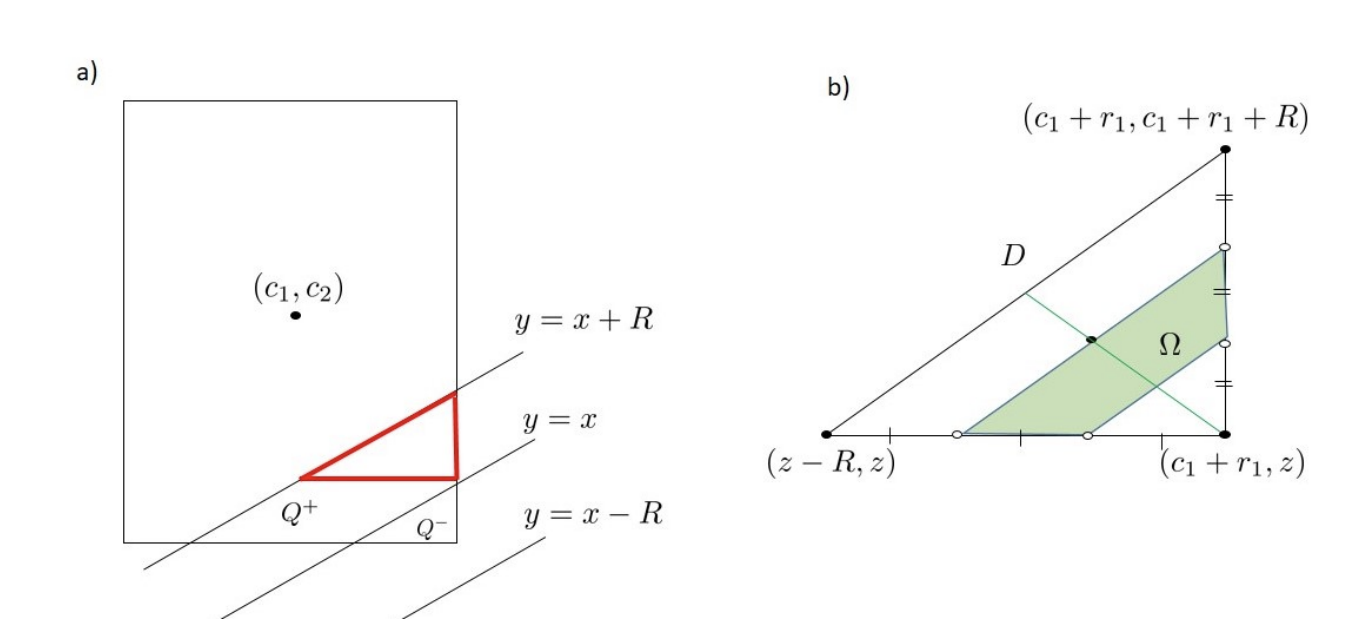}
\end{center}
Figure 2: a) Illustration of the rectangle $Q$. b) Illustration of the trapezium $\Omega$ contained in D.\\[3mm]
Now, denoting

$$
\begin{aligned}
\varphi^{*}\left(c_{1}, r_{1}, c_{2}, r_{2}, R\right) & :=c_{1}+r_{1}-z+R \\
& =\min \left\{r_{1}-r_{2}+c_{2}-c_{1}+R, c_{1}-c_{2}+r_{1}+r_{2}+R, R\right\}
\end{aligned}
$$

by our assumptions we have that $\varphi^{*}\left(c_{1}, r_{1}, c_{2}, r_{2}, R\right)>0$. Hence, we may conclude that

\begin{equation*}
\frac{d^{+}}{d \tau} I(0)>\frac{1}{6} \varphi^{*}\left(c_{1}, r_{1}, c_{2}, r_{2}, R\right) \min _{r \in\left[\frac{R}{3 \sqrt{2}}, \frac{R}{\sqrt{2}}\right]} K^{\prime}(r)>0 \tag{16}
\end{equation*}

For $r_{1}>r_{2}$, one can obtain in a similar way the inequality (16), where $r_{1}$ and $r_{2}$ are swapped. Accommodating also the case of $c_{1}>c_{2}$ the inequality (16) holds with $\varphi^{*}\left(c_{1}, r_{1}, c_{2}, r_{2}, R\right)=\varphi\left(c_{1}, r_{1}, c_{2}, r_{2}, R\right)$, as defined in (14)

The proof of Lemma 2.12 below follows the same as that of [17, Lemma 2.17 (a)] and will therefore be omitted.

Lemma 2.12. Assume $K \in C^{1}(\mathbb{R})$ is as defined in Lemma 1. For any open sets $U_{1}, U_{2} \subset \mathbb{R}$, let $\mu_{i}:=\chi_{U_{i}}$ for $i=1,2$ and $I(\tau)=I_{K}\left[\mu_{1}, \mu_{2}\right](\tau)$. Then

$$
\frac{d}{d \tau} I(\tau) \geq 0 \text { for all } \tau \geq 0
$$

Proof of Theorem 2.10.1.\\
Assume, for a contradiction, that $\left.\mu_{0}\right|_{D}$ is not radially decreasing with respect to any $x_{0} \in D$ considered as a centre. Then by [17, Lemma 2.18] there exists a unit vector $e$ such that $\mu_{0}$ is not symmetric decreasing about any hyperplane with normal vector $e$. We set $e=(1,0, \ldots, 0)$ without loss of generality.

Recall Definition 2.6, where we define the continuous Steiner symmetrization of a non-negative function on $\mathbb{R}^{d}$ with respect to a specific direction. In order to prove Theorem 2.10.1, we modify $S^{\tau} \mu_{0}$ in such a way that $U_{x^{\prime}}^{h}$ travels at the speed $v(h)$ instead of a constant speed 1 , where

$$
v(h):= \begin{cases}1, & \text { if } h \geq h_{0} \\ \left(\frac{h}{h_{0}}\right)^{m-1}, & \text { if } 0<h<h_{0}\end{cases}
$$

for some $h_{0}$ sufficiently small, defined later. We let $\mu(\tau, \cdot)=\tilde{S}^{\tau} \mu_{0}$ where

$$
\tilde{S}^{\tau} \mu_{0}:=\int_{0}^{\infty} \chi_{M^{v(h) \tau}\left(U_{x^{\prime}}^{h}\right)}\left(x_{1}\right) d h
$$

Now, in [17, Proposition 2.8] it is shown that there exists some $\delta_{1}>0$ and $C>0$, depending on $m, \mu_{0}$, and $W$, such that, for any $\tau \in\left[0, \delta_{1}\right]$,

\begin{align*}
\left|\mu(\tau, x)-\mu_{0}(x)\right| & \leq C \mu_{0}(x) \tau \text { for all } x \in \mathbb{R}^{d}  \tag{17}\\
\int_{D_{i}}\left(\mu(\tau, x)-\mu_{0}(x)\right) d x & =0 \tag{18}
\end{align*}

for any connected component $D_{i}$ of $\operatorname{supp} \mu_{0}$. The proof of [17, Proposition 2.8] does not rely on the infinite support of $\omega^{\prime}$ assumed in [17]. As a result, we see that (17) and (18) hold for $W$ satisfying assumptions ( $W 1$ ), (W2), and (W4). Hence, we can use the same argument as in the proof of [17, Theorem 2.2] to obtain that there exists some $C_{2}>0$ and $\delta_{0}>0$ with $\delta_{0} \leq \delta_{1}$ such that

\begin{equation*}
\left|\mathcal{E}[\mu(\tau, \cdot)]-\mathcal{E}\left[\mu_{0}\right]\right| \leq C_{2} \tau^{2} \text { for all } \tau \in\left[0, \delta_{0}\right] \tag{19}
\end{equation*}

It remains to be shown that there exists a $C_{1}>0$ and some $\tau_{1}, \tau_{2}$ with $0 \leq \tau_{1}<\tau_{2} \leq \delta_{0}$ such that

\begin{equation*}
\mathcal{E}[\mu(\tau, \cdot)]-\mathcal{E}\left[\mu_{0}\right] \leq-C_{1} \tau \text { for all } \tau \in\left[\tau_{1}, \tau_{2}\right] \tag{20}
\end{equation*}

This will allow us to conclude that (19) and (20) hold for all $\tau \in\left[\tau_{1}, \tau_{2}\right]$. Combining (19) and (20) will then lead to a contradiction of $\tau_{1}$.

In the proof of $[17$, Proposition 2.8$]$, it is shown that $\mathcal{S}[\mu(\tau, \cdot)] \leq \mathcal{S}\left[\mu_{0}\right]$ for all $\tau>0$. Hence, it is sufficient to show that

$$
\mathcal{I}[\mu(\tau, \cdot)]-\mathcal{I}\left[\mu_{0}\right] \leq-C_{1} \tau \text { for all } \tau \in\left[\tau_{1}, \tau_{2}\right]
$$

We prove this result as follows: Fix $\tau \in\left[\tau_{1}, \tau_{2}\right]$, where $\tau_{1}, \tau_{2}$ are to be defined later. Fix $x^{\prime} \in D_{x_{1}}$, where $D_{x_{1}}=\left\{z \in \mathbb{R}^{d-1}:\left(x_{1}, z\right) \in D\right\}$. Consider the interval $[a, b]$, where

$$
a:=\max \left\{x_{1} \in \mathbb{R}: \mu_{0}\left(x_{1}, x^{\prime}\right)=\max _{x \in \mathbb{R}} \mu_{0}\left(x, x^{\prime}\right)\right\}
$$

and

$$
b:=\min \left\{x_{1} \in\left[a, a+\frac{R}{2}\right]: \mu_{0}\left(x_{1}, x^{\prime}\right)=\min _{x \in\left[a, a+\frac{R}{2}\right]} \mu_{0}\left(x, x^{\prime}\right)\right\}
$$

From the definition of $a$ and $b$, we have that $0<b-a \leq \frac{R}{2}$ and $\mu_{0}\left(a, x^{\prime}\right)>$ $\mu_{0}\left(x, x^{\prime}\right)>\mu_{0}\left(b, x^{\prime}\right)$ for all $x \in(a, b)$. For $\alpha>0$, define $\mathcal{H}^{\alpha}:=\left\{\left(x^{\prime}, h\right) \in\right.$ $\left.D_{x_{1}} \times(0, \infty):\left|U_{x^{\prime}}^{h} \cap[a, b]\right|>\alpha\right\}$. By continuity of $\mu_{0}$, we can choose $\alpha$ sufficiently small so that $\mathcal{H}^{\alpha}$ has positive measure. By [17, Proposition 2.8], for any $h_{0}>0$, $\mu(\tau, \cdot)$ satisfies (17) and (18) for all $\tau \in\left[\tau_{1}, \tau_{2}\right]$, since $\tau_{1}, \tau_{2} \leq \delta_{1}$. Hence, we may choose $h_{0}:=\mu_{0}\left(\frac{a+b}{2}, x^{\prime}\right)$ and define

$$
B_{1}=\left\{\left(x^{\prime}, h\right) \in \mathcal{H}^{\alpha}: h \geq h_{0}\right\} \text { and } B_{2}=\left\{\left(x^{\prime}, h\right) \in \mathcal{H}^{\alpha}: h<h_{0}\right\}
$$

By our choice of $h_{0}$ and continuity of $\mu_{0}$ we have that $B_{1}$ and $B_{2}$ have positive measure. Now, consider $\left(x^{\prime}, h_{1}\right) \in B_{1}$ and $\left(y^{\prime}, h_{2}\right) \in B_{2}$ where $0<\left\|x^{\prime}-y^{\prime}\right\|<1$. Let $U_{x^{\prime}}^{h_{1}}$ and $U_{y^{\prime}}^{h_{2}}$ consist of either a finite or infinite union of disjoint open intervals. Then, there exist intervals $I\left(c_{k *}^{1}, r_{k *}^{1}\right)$ and $I\left(c_{\ell *}^{2}, r_{\ell *}^{2}\right)$ in $U_{x^{\prime}}^{h_{1}}$ and $U_{y^{\prime}}^{h_{2}}$, respectively, such that $\left|I\left(c_{k *}^{1}, r_{k *}^{1}\right) \cap[a, b]\right|>\beta$ and $\left|I\left(c_{\ell *}^{2}, r_{\ell *}^{2}\right) \cap[a, b]\right|>\beta$ for some $\beta \in(0, \alpha]$.

We show that there exists some $\tau_{1}, \tau_{2}$ with $0 \leq \tau_{1}<\tau_{2} \leq \delta_{0}$ such that $\mu_{1}:=\chi_{M^{v\left(h_{1}\right)} \tau_{0}\left(I\left(c_{k *}^{1}, r_{k *}^{1}\right)\right)}$ and $\mu_{2}:=\chi_{M^{v\left(h_{2}\right) \tau_{0}\left(I\left(c_{\ell *}^{2}, r_{\ell *}^{2}\right)\right)}}$ satisfy the assumptions of Lemma 2.11.2 for all $\tau_{0} \in\left[\tau_{1}, \tau_{2}\right]$. We choose $H=\left\{x_{1}=\tilde{x}\right\}$ for some $\tilde{x} \in D_{x^{\prime}}=\left\{x_{1} \in \mathbb{R}:\left(x_{1}, x^{\prime}\right) \in D\right\}$ and let $I\left(c_{k^{*}}^{1}, r_{k^{*}}^{1}\right)=: I\left(\tilde{x}+c_{1}, r_{1}\right)$ and $I\left(c_{k *}^{2}, r_{k^{*}}^{2}\right)=: I\left(\tilde{x}+c_{2}, r_{2}\right)$. We consider the following two cases: $c_{1}=c_{2}$ and $c_{1} \neq c_{2}$.

Case 1: $c_{1}=c_{2}$\\
We choose $\tilde{x}$ such that $\left|c_{1}\right|=\left|\tilde{x}-c_{k *}^{1}\right|=\frac{1}{4} \min \left\{R, \delta_{0}, \tau^{*}, \frac{C_{1}}{C_{2}}\right\}$ where $\tau^{*}$ is the value at which $M^{v\left(h_{2}\right) \tau^{*}}\left(I\left(c_{k *}^{2}, r_{k^{*}}^{2}\right)\right)$ shares a common endpoint with a neighbouring interval (if one exists). Since $\mu_{0}$ is not symmetric decreasing about any hyperplane with normal vector $e$, we have that $\mu_{0}$ is not symmetric decreasing about $H=\left\{x_{1}=\tilde{x}\right\}$. Now, let $\tau_{1}=\left|c_{1}\right|$. We claim that the assumptions of Lemma 2.11.2 hold at $\tau=\tau_{1}$. Indeed, we have that

$$
M^{v\left(h_{1}\right) \tau_{1}}\left(I\left(\tilde{x}+c_{1}, r_{1}\right)\right)=I\left(\tilde{x}+c_{1}-v\left(h_{1}\right) \tau_{1} \operatorname{sgn} c_{1}, r_{1}\right)
$$

and

$$
M^{v\left(h_{2}\right) \tau_{1}}\left(I\left(\tilde{x}+c_{2}, r_{2}\right)\right)=I\left(\tilde{x}+c_{2}-v\left(h_{2}\right) \tau_{1} \operatorname{sgn} c_{2}, r_{2}\right)
$$

Set $\hat{c}_{1}:=c_{1}-v\left(h_{1}\right) \tau_{1} \operatorname{sgn} c_{1}$ and $\hat{c}_{2}:=c_{2}-v\left(h_{2}\right) \tau_{1} \operatorname{sgn} c_{2}$. Then, $\hat{c}_{1}=0$, since $v\left(h_{1}\right)=1$. If $c_{2}>0$, then $\hat{c}_{2}=\left|c_{2}\right|-v\left(h_{2}\right)\left|c_{2}\right|>0$ and if $c_{2}<0$, then $\hat{c}_{2}=-\left|c_{2}\right|+v\left(h_{2}\right)\left|c_{2}\right|<0$. This follows from the fact that $v\left(h_{2}\right)=\left(\frac{h_{2}}{h_{0}}\right)^{m-1}<$ 1. Hence, $\operatorname{sgn} \hat{c}_{1} \neq \operatorname{sgn} \hat{c}_{2}$. Furthermore, since $h_{1}>h_{2}$ implies $r_{1}<r_{2}$, by continuity of $\mu_{0}$, and since $\left|c_{1}\right| \leq \frac{R}{4}$, we have that\\
$\left|\hat{c}_{2}-\hat{c}_{1}\right|=\left|\hat{c}_{2}\right|=\left|c_{2}-v\left(h_{2}\right) \tau_{1} \operatorname{sgn} c_{2}\right| \leq\left|c_{2}\right|+v\left(h_{2}\right)\left|c_{1}\right|<2\left|c_{1}\right| \leq \frac{R}{2}<r_{1}+r_{2}+R$.\\[3mm]
\includegraphics[max width=\textwidth, center]{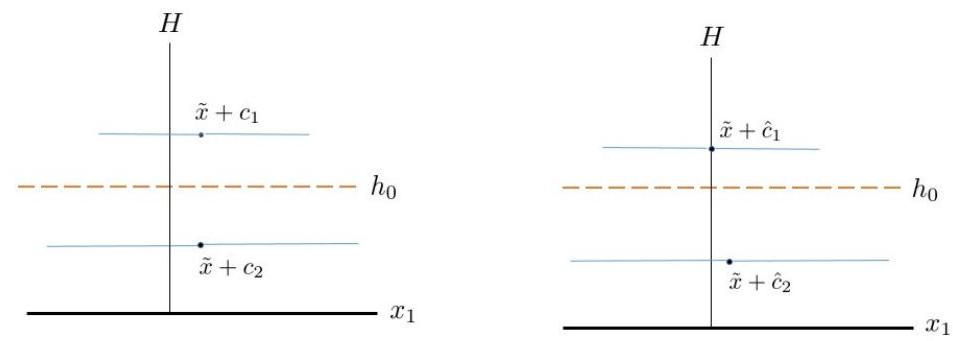}\\
Figure 3: Illustration of case 1 where $c_{1}=c_{2}>0$ and $x^{\prime}$ is fixed. On the left, $\tau=0$ and on the right, $\tau=\tau_{1}$.\\[3mm]
Also, since $\left[\tilde{x}+c_{1}+r_{1}, \tilde{x}+c_{2}+r_{2}\right] \subseteq[a, b]$ and $b-a \leq \frac{R}{2}$, it follows that

$$
r_{2}-r_{1}=\left(\tilde{x}+c_{2}+r_{2}\right)-\left(\tilde{x}+c_{1}+r_{1}\right) \leq \frac{R}{2}<\left|\hat{c}_{2}-\hat{c}_{1}\right|+R
$$

as required.\\
It remains to be shown that the assumptions of Lemma 2.11 .2 continue to hold for all $\tau \in\left[\tau_{1}, \tau_{2}\right]$, where $\tau_{2}:=\min \left\{\tau_{1}+\left|\hat{c}_{2}\right|, \frac{\tau^{*}}{2}, \delta_{0}\right\}$. Now, $M^{v\left(h_{1}\right) \tau_{2}}(I(\tilde{x}+$ $\left.\left.c_{1}, r_{1}\right)\right)=I\left(\tilde{x}, r_{1}\right)$ and $M^{v\left(h_{2}\right) \tau_{2}}\left(I\left(\tilde{x}+c_{2}, r_{2}\right)\right)=I\left(\tilde{x}+c_{2}-v\left(h_{2}\right) \tau_{2} \operatorname{sgn} c_{2}, r_{2}\right)$.

Set $\tilde{c}_{1}:=c_{1}-v\left(h_{1}\right) \tau_{2} \operatorname{sgn} c_{1}$ and $\tilde{c}_{2}:=c_{2}-v\left(h_{2}\right) \tau_{2} \operatorname{sgn} c_{2}$. Then $\tilde{c}_{1}=0$. Furthermore, if $\hat{c}_{2}>0$ then $c_{2}>0$ and so $\tilde{c}_{2}=c_{2}-v\left(h_{2}\right) \tau_{2} \geq\left|\hat{c}_{2}\right|-v\left(h_{2}\right)\left|\hat{c}_{2}\right|>0$ and if $\hat{c}_{2}<0$ then $c_{2}<0$ and so $\tilde{c}_{2}=c_{2}+v\left(h_{2}\right) \tau_{2} \leq-\left|\hat{c}_{2}\right|+v\left(h_{2}\right)\left|\hat{c}_{2}\right|<0$. Hence, $\operatorname{sgn} \tilde{c}_{1} \neq \operatorname{sgn} \tilde{c}_{2}$. Additionally, since $\left|\hat{c}_{2}\right| \leq \frac{R}{2}$ by (21), we have that

$$
\left|\tilde{c}_{2}-\tilde{c}_{1}\right|=\left|\tilde{c}_{2}\right| \leq\left|c_{2}\right|+v\left(h_{2}\right) \tau_{2} \leq\left|c_{2}\right|+v\left(h_{2}\right)\left(\left|c_{2}\right|+\left|\hat{c}_{2}\right|\right)<r_{2}+r_{1}+R
$$

and

$$
\left|r_{2}-r_{1}\right| \leq \frac{R}{2}<\left|\tilde{c}_{2}-\tilde{c}_{1}\right|+R
$$

as required.\\
Case 2: $c_{1} \neq c_{2}$\\
Suppose $c_{1} \neq c_{2}$. We choose $\tilde{x}=c_{k^{*}}^{1}$ such that $c_{1}=0$ and $\operatorname{sgn} c_{1} \neq \operatorname{sgn} c_{2}$. Now, let $\tau_{1}=0$ and $\tau_{2}=\frac{1}{2} \min \left\{\left|c_{2}\right|, \tau^{*}, \delta_{0}, R\right\}>0$. Then we have that $M^{v\left(h_{1}\right) \tau_{2}}\left(I\left(\tilde{x}+c_{1}, r_{1}\right)\right)=I\left(\tilde{x}+\tilde{c}_{1}, r_{1}\right)$ and $M^{v\left(h_{2}\right) \tau_{2}}\left(I\left(\tilde{x}+c_{2}, r_{2}\right)\right)=I\left(\tilde{x}+\tilde{c}_{2}, r_{2}\right)$ where $\tilde{c}_{1}=0$ and $\tilde{c}_{2}=c_{2}-v\left(h_{2}\right) \tau_{2} \operatorname{sgn} c_{2}$. If $c_{2}>0$, then $\tilde{c}_{2}=c_{2}-v\left(h_{2}\right) \tau_{2} \geq$ $\left|c_{2}\right|-v\left(h_{2}\right)\left|c_{2}\right|>0$. Also, if $c_{2}<0$, then $\tilde{c}_{2}=c_{2}+v\left(h_{2}\right) \tau_{2} \leq-\left|c_{2}\right|+v\left(h_{2}\right)\left|c_{2}\right|<0$.

Hence, $\operatorname{sgn} \tilde{c}_{1} \neq \operatorname{sgn} \tilde{c}_{2}$. Furthermore, by continuity of $\mu_{0}$, we have that $a \leq \tilde{x}+c_{1}+r_{1}=\tilde{x}+\tilde{c}_{1}+r_{1}<\tilde{x}+c_{2}+r_{2} \leq b$. Now, if $c_{2}>0$, then $a \leq \tilde{x}+\tilde{c}_{1}+r_{1}=\tilde{x}+r_{1}<\tilde{x}+\tilde{c}_{2}+r_{1}<\tilde{x}+\tilde{c}_{2}+r_{2}<b$. Hence, since $b-a \leq \frac{R}{2}$,\\
we have that $\tilde{x}+\tilde{c}_{2}+r_{2}-\left(\tilde{x}+\tilde{c}_{1}+r_{1}\right) \leq \frac{R}{2}$ and so $\tilde{c}_{2}-\tilde{c}_{1} \leq r_{1}-r_{2}+\frac{R}{2}<r_{1}+r_{2}+R$.\\
If $c_{2}<0$, then $a \leq \tilde{x}+\tilde{c}_{1}+r_{1}=\tilde{x}+r_{1}<\tilde{x}+\tilde{c}_{2}+r_{2}=\tilde{x}+c_{2}+v\left(h_{2}\right) \tau_{2}+r_{2} \leq$ $\tilde{x}+c_{2}+v\left(h_{2}\right) \frac{R}{2}+r_{2} \leq b+v\left(h_{2}\right) \frac{R}{2}$. Hence, $b+v\left(h_{2}\right) \frac{R}{2}-a \leq \frac{R}{2}+v\left(h_{2}\right) \frac{R}{2}<R$ and so $\tilde{x}+\tilde{c}_{2}+r_{2}-\left(\tilde{x}+\tilde{c}_{1}+r_{1}\right)<R$, which yields $\tilde{c}_{2}-\tilde{c}_{1}<r_{1}+r_{2}+R$.

Furthermore, for both $c_{2}>0$ and $c_{2}<0$, we see that

$$
r_{2}-r_{1}<\tilde{c}_{1}-\tilde{c}_{2}+R \leq\left|\tilde{c}_{2}-\tilde{c}_{1}\right|+R
$$

as required. Hence, for both cases 1 and 2 it follows that $\chi_{M^{v\left(h_{1}\right) \tau_{0}\left(I\left(c_{k^{*}}^{1}, r_{k^{*}}^{1}\right)\right)}}$ and $\chi_{M^{v\left(h_{2}\right) \tau_{0}\left(I\left(c_{k^{*}}^{2}, r_{k^{*}}^{2}\right)\right)}}$ satisfy the asumptions of Lemma 2.11.2 for all $\tau_{0} \in\left[\tau_{1}, \tau_{2}\right]$. Thus, for $s$ sufficiently small and $\tau=v(h) \tau_{0}+s$, we have that

\begin{center}
\includegraphics[max width=\textwidth]{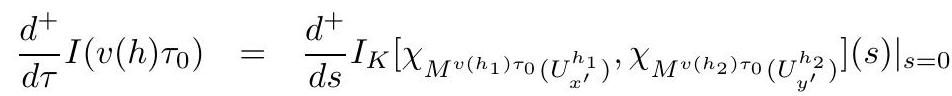}
\end{center}

\begin{center}
\includegraphics[max width=\textwidth]{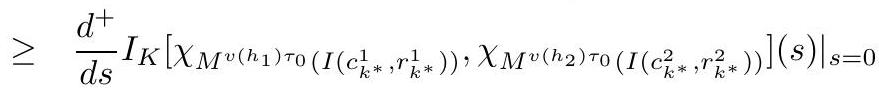}
\end{center}

\begin{align*}
& \geq \quad c>0, \tag{22}
\end{align*}

for all $\tau_{0} \in\left[\tau_{1}, \tau_{2}\right]$, by Lemma 2.11 and 2.12 .\\
Now, define a 1-D kernel $K_{l}(z)=-\frac{1}{2} \omega\left(\sqrt{z^{2}+l^{2}}\right)$. For any $l>0, K_{l} \in C^{1}(\mathbb{R})$ is even with $K_{l}(z)<0$ for all $0<z<R$, where $R:=\sqrt{1-l^{2}}, l^{2}<1$, and $K_{l}(z)=0$ for all $z \geq R$. Then,

$$
\begin{aligned}
\mathcal{I}\left[\tilde{S}^{\tau} \mu_{0}\right] & =\frac{1}{2} \int_{D} \int_{D} \tilde{S}^{\tau} \mu_{0}(x) \tilde{S}^{\tau} \mu_{0}(y) W(x-y) d y d x \\
& =\frac{1}{2} \int_{D^{2}} \int_{\left(\mathbb{R}^{+}\right)^{2}} \chi_{M^{v\left(h_{1}\right) \tau}\left(U_{x^{\prime}}^{h_{1}}\right)}\left(x_{1}\right) \chi_{M^{v\left(h_{2}\right) \tau}\left(U_{y^{\prime}}^{h_{2}}\right)}\left(y_{1}\right) W(x-y) d h_{1} d h_{2} d y d x \\
& =-\int_{D^{2}} \int_{\left(\mathbb{R}^{+}\right)^{2}} \chi_{M^{v\left(h_{1}\right) \tau}\left(U_{x^{\prime}}^{h_{1}}\right)}\left(x_{1}\right) \chi_{M^{v\left(h_{2}\right) \tau}\left(U_{y^{\prime}}^{h_{2}}\right)}\left(y_{1}\right) K_{\| x^{\prime}-y^{\prime}||}\left(\left|x_{1}-y_{1}\right|\right) d h_{1} d h_{2} d y d x
\end{aligned}
$$

This follows from the definition of $\tilde{S}^{\tau} \mu_{0}$ and since $W(x-y)=\omega(\|x-y\|)=$ $\omega\left(\left(\left(x_{1}-y_{1}\right)^{2}+\ldots+\left(x_{d}-y_{d}\right)^{2}\right)^{1 / 2}\right)=-2 K_{\left\|x^{\prime}-y^{\prime}\right\|}\left(\left|x_{1}-y_{1}\right|\right)$, where $0<\left\|x^{\prime}-y^{\prime}\right\|<$ 1. Now, using the definition given in (12), we have that

$$
\mathcal{I}\left[\tilde{S}^{\tau} \mu_{0}\right]=-\int_{D_{x_{1}}} \int_{D_{y_{1}}} \int_{\left(\mathbb{R}^{+}\right)^{2}} I_{K_{\left\|x^{\prime}-y^{\prime}\right\|}}\left[\chi_{U_{x^{\prime}}^{h_{1}}}, \chi_{U_{y^{\prime}}^{h_{2}}}\right](v(h) \tau) d h_{1} d h_{2} d y^{\prime} d x^{\prime}
$$

Taking the right derivative, we obtain

$$
\begin{aligned}
-\frac{d^{+}}{d \tau} \mathcal{I}\left[\tilde{S}^{\tau} \mu_{0}\right] & =\int_{D_{x_{1}}} \int_{D_{y_{1}}} \int_{\left(\mathbb{R}^{+}\right)^{2}} \frac{d^{+}}{d \tau} I_{K_{\left\|x^{\prime}-y^{\prime}\right\|}}\left[\chi_{U_{x^{\prime}}^{h_{1}}}, \chi_{\left.U_{y^{\prime}}^{h_{2}}\right]}\right](v(h) \tau) d h_{1} d h_{2} d y^{\prime} d x^{\prime} \\
& \geq \int_{B_{1}} \int_{B_{2}} \frac{d^{+}}{d \tau} I_{K_{\left\|x^{\prime}-y^{\prime}\right\| \mid}}\left[\chi_{U_{x^{\prime}}^{h_{1}}}, \chi_{U_{y^{\prime}}^{h_{2}}}\right](v(h) \tau) d y^{\prime} d h_{2} d x^{\prime} d h_{1}
\end{aligned}
$$

by Lemma 2.12. Now, for any $\left(x^{\prime}, h_{1}\right) \in B_{1}$ and $\left(y^{\prime}, h_{2}\right) \in B_{2}$, by (22), we have that

$$
-\frac{d^{+}}{d \tau} \mathcal{I}\left[\tilde{S}^{\tau} \mu_{0}\right] \geq\left|B_{1}\right|\left|B_{2}\right| c=: C_{1}>0
$$

for all $\tau \in\left[\tau_{1}, \tau_{2}\right]$. Hence, by the fundamental theorem of calculus and using the fact that $\tilde{S}^{0} \mu_{0}=\mu_{0}$, it follows that $\mathcal{I}\left[\tilde{S}^{\tau} \mu_{0}\right]-\mathcal{I}\left[\mu_{0}\right] \leq-C_{1} \tau$, for all $\tau \in\left[\tau_{1}, \tau_{2}\right]$. Hence, we have that

$$
\mathcal{E}[\mu(\tau, \cdot)]-\mathcal{E}\left[\mu_{0}\right] \leq-C_{1} \tau \quad \text { and }\left|\mathcal{E}[\mu(\tau, \cdot)]-\mathcal{E}\left[\mu_{0}\right]\right| \leq C_{2} \tau^{2}
$$
for all $\tau \in\left[\tau_{1}, \tau_{2}\right]$ (since $\tau_{2} \leq \delta_{0}$ ). Now, by definition, $\tau_{1}<\frac{C_{1}}{2 C_{2}}$. But, since $-C_{2} \tau_{1}^{2} \leq \mathcal{E}[\mu(\tau, \cdot)]-\mathcal{E}\left[\mu_{0}\right] \leq-C_{1} \tau_{1}$, it must be that $\tau_{1} \geq \frac{C_{1}}{C_{2}}$, a contradiction.

We show in Theorem 2.13 below that the compactness of each connected component of $\operatorname{supp} \rho_{s}$ follows from Theorem 2.10.1.

Theorem 2.13. Let $m>2$. Assume ( $W 1$ ), (W2), and (W4) hold. If $\rho_{s}$ is a stationary state of (1), then each connected component of $\operatorname{supp} \rho_{s}$ is compact.

Proof. We consider an arbitrary connected component $D$ of $\operatorname{supp} \rho_{s}$. We assume, for a contradiction, that $D$ is not compact. Since $D$ is closed, it must be that $D$ is unbounded. Therefore, since $\rho_{s}$ is radially symmetric we have that $D=\mathbb{R}^{d}$. The proof then follows exactly the same as that of Theorem 2.5.

Since each connected component $D$ of the support of $\rho_{s}$ is compact by Theorem 2.13, we have that $\left.\rho_{s}\right|_{\partial D}=0$. Hence, as in the statement of Theorem 2.10.2, we consider $x \in \operatorname{interior}(D)$ and $y \in \operatorname{interior}\left(\operatorname{supp} \rho_{s} \backslash D\right)$. Note that, since $D$ and closure $\left(\operatorname{supp} \rho_{s} \backslash D\right)$ are disjoint by assumption and are compact, we have that $\min \left\{\|x-y\|: x \in D, y \in \operatorname{closure}\left(\operatorname{supp} \rho_{s} \backslash D\right)\right\}>0$.

Proof of Theorem 2.10.2. Set $\rho_{s}=\mu_{0}=\mu(0, \cdot)$. We consider supp $\mu_{0}$ made up of two connected components, $D_{1}$ and $D_{2}$. Suppose for a contradiction that there exists some $x^{*}=\left(x_{1}^{*}, x^{\prime *}\right) \in \operatorname{interior} D_{1}$ and $y^{*}=\left(y_{1}^{*}, y^{\prime *}\right) \in$ interior $D_{2}$ such that $\left\|x^{*}-y^{*}\right\|<1$. Then, $\left|x_{1}^{*}-y_{1}^{*}\right|<\sqrt{1-\left\|x^{\prime *}-y^{\prime *}\right\|^{2}}=R$, where $\left\|x^{* *}-y^{\prime *}\right\|^{2}<1$

Set $x_{1}^{*}=\tilde{x}+c_{1}+r_{1}$ and $y_{1}^{*}=\tilde{x}+c_{2}-r_{2}$ where $\tilde{x}$ is chosen such that $c_{1}<0$ and $c_{2}>0$. We claim that $\chi_{I\left(\tilde{x}+c_{1}, r_{1}\right)}$ and $\chi_{I\left(\tilde{x}+c_{2}, r_{2}\right)}$ satisfy the assumptions of Lemma 2.11.2. This can be seen from the fact that $\left|x_{1}^{*}-y_{1}^{*}\right|<R$ implies $\left|c_{2}-r_{2}-c_{1}-r_{1}\right|<R$, which gives ||$c_{2}-c_{1}|-| r_{2}+r_{1}||<R$, yielding $\left|c_{2}-c_{1}\right|<r_{1}+r_{2}+R$ and $\left|r_{2}-r_{1}\right|<\left|c_{2}-c_{1}\right|+R$.

Furthemore, taking $\tau_{m}:=\frac{1}{2} \min \left\{\left|c_{1}\right|,\left|c_{2}\right|,\left|x_{1}^{*}-y_{1}^{*}\right|, \delta_{0}\right\}>0$, with $\delta_{0}$ as defined in [17, Theorem 2.2], we have that $M^{\tau_{m}}\left(I\left(\tilde{x}+c_{1}, r_{1}\right)\right)=I\left(\tilde{x}+\tilde{c}_{1}, r_{1}\right)$\\
and $M^{\tau_{m}}\left(I\left(\tilde{x}+c_{2}, r_{2}\right)\right)=I\left(\tilde{x}+\tilde{c}_{2}, r_{2}\right)$ where $\tilde{c}_{1}=c_{1}-\tau_{m} \operatorname{sgn} c_{1}=c_{1}+\tau_{m}$ $<-\left|c_{1}\right|+\left|c_{1}\right|=0$, and $\tilde{c}_{2}=c_{2}-\tau_{m} \operatorname{sgn} c_{2}=c_{2}-\tau_{m}>\left|c_{2}\right|-\left|c_{2}\right|=0$. Hence, $\operatorname{sgn} \tilde{c}_{1} \neq \operatorname{sgn} \tilde{c}_{2}$

Also, since $\left[\tilde{x}+\tilde{c}_{1}+r_{1}, \tilde{x}+\tilde{c}_{2}-r_{2}\right] \subset\left[\tilde{x}+c_{1}+r_{1}, \tilde{x}+c_{2}-r_{2}\right]$, we have that $\left|\tilde{c}_{2}-\tilde{c}_{1}-r_{2}-r_{1}\right|<\left|c_{2}-c_{1}-r_{2}-r_{1}\right|<R$, yielding $\left|\tilde{c}_{2}-\tilde{c}_{1}\right|<r_{1}+r_{2}+R$ and $\left|r_{2}-r_{1}\right|<\left|\tilde{c}_{2}-\tilde{c}_{1}\right|+R$, as required.

Furthermore, setting $h_{0}=\max \mu_{0}$, we have that $\chi_{M^{v(h) \tau_{m}\left(I\left(\tilde{x}+c_{1}, r_{1}\right)\right)}}$ and $\chi_{M^{v(h) \tau_{m}}\left(I\left(\tilde{x}+c_{2}, r_{2}\right)\right)}$ satisfy the assumptions of Lemma 2.11.2. for all $\tau \in\left[0, \tau_{m}\right]$, since $v(h)<1$ for $h<h_{0}$.

Now, fix $\left(x^{\prime}, h_{1}\right) \in D_{1}^{x_{1}} \times(0, \infty)$ and $\left(y^{\prime}, h_{2}\right) \in D_{2}^{y_{1}} \times(0, \infty)$, where $D_{i}^{z_{1}}=$ $\left\{z^{\prime} \in \mathbb{R}^{d-1}:\left(z_{1}, z^{\prime}\right) \in D_{i}\right\}$, for $i=1,2$ and $z_{1} \in \mathbb{R}$. We define the set

$$
B:=\left\{\left(x^{\prime}, h_{1}\right) \times\left(y^{\prime}, h_{2}\right):|\alpha-\beta|<R, \alpha \in \operatorname{interior} D_{1}^{x^{\prime}}, \beta \in \operatorname{interior} D_{2}^{y^{\prime}}\right\}
$$

where $D_{i}^{z^{\prime}}=\left\{z_{1} \in \mathbb{R}:\left(z_{1}, z^{\prime}\right) \in D_{i}\right\}$, for $i=1,2$ and $z^{\prime} \in \mathbb{R}^{d-1}$. We know by assumption that $B$ is non-empty. Furthermore, since $B$ is an open set, it has positive measure. Hence, using the same approach as in the proof of Theorem 2.10.1, we find that

$$
\begin{aligned}
-\frac{d^{+}}{d \tau} \mathcal{I}\left[\tilde{S}^{\tau} \mu_{0}\right] & =\int_{\left(\mathbb{R}^{d-1}\right)^{2}} \int_{\left(\mathbb{R}^{+}\right)^{2}} \frac{d^{+}}{d \tau} I_{K_{\left\|x^{\prime}-y^{\prime}\right\|}}\left[\chi_{U_{x^{\prime}}^{h_{1}}}, \chi_{U_{y^{\prime}}^{h_{2}}}\right](v(h) \tau) d h_{1} d h_{2} d x^{\prime} d y^{\prime} \\
& \geq c|B|>0 \text { for all } \tau \in\left[0, \tau_{m}\right]
\end{aligned}
$$

by Lemma 2.11 and 2.12. Thus, there exists a $C>0$ such that

\begin{equation*}
\mathcal{E}\left[\tilde{S}^{\tau} \mu_{0}\right]-\mathcal{E}\left[\mu_{0}\right]<-C \tau<0 \tag{23}
\end{equation*}

for all $\tau \in\left(0, \tau_{m}\right]$. However; using the fact that $\mu_{0}$ is a stationary state of (1) we have that

\begin{equation*}
\left|\mathcal{E}[\mu(\tau, \cdot)]-\mathcal{E}\left[\mu_{0}\right]\right| \leq C_{2} \tau^{2} \tag{24}
\end{equation*}

for all $\tau \in\left[0, \tau_{m}\right]$, as discussed in the proof of Theorem 2.10.1. Combining Inequality (23) and (24) and taking $\tau=\frac{C}{2 C_{2}}$, we see that
$$
-\frac{C^{2}}{4 C_{2}} \leq \mathcal{E}[\mu(\tau, \cdot)]-\mathcal{E}\left[\mu_{0}\right]<-\frac{C^{2}}{2 C_{2}}
$$
a contradiction. Hence, from the theorems in this section, we have provided a proof for Theorem 2.1.

\section*{3. Mass-independent boundedness of stationary states}
In this section, we prove that for $m>2$ and for attractive kernels with both compact and infinite support, stationary states of Equation (1) possess a mass-independent upper-bound. That is, regardless of the size of the support of $W$, there exists an upper-bound of any stationary solution to (1) that does not depend on the initial condition.

Theorem 3.1. Let $m>2$. If $\rho_{s}$ is a stationary state of (1) with $W$ satisfying assumptions (W1), (W2), and (W4), then
$$
\rho_{s}(x) \leq \rho_{s}^{*}:=\left(\frac{m-1}{m \varepsilon}\|W\|_{L^{1}\left(\mathbb{R}^{d}\right)}\right)^{1 /(m-2)} \text { for all } x \in \mathbb{R}^{d}
$$

Proof. Since $\rho_{s}$ is a stationary state of (1),

$$
\frac{m \varepsilon}{m-1} \rho_{s}^{m-1}(x)+\left(W * \rho_{s}\right)(x)=C_{j} \text { for all } x \in D_{j}
$$
where $C_{j}$ may be different on each connected component $D_{j}$ of $\operatorname{supp} \rho_{s}$. Now, let $z \in \partial D_{j}$. Then,
$$
\frac{m \varepsilon}{m-1} \rho_{s}^{m-1}(x)+\left(W * \rho_{s}\right)(x)=\frac{m \varepsilon}{m-1} \rho_{s}^{m-1}(z)+\left(W * \rho_{s}\right)(z) \text { for all } x \in D_{j}
$$
Since $\rho_{s} \in C\left(\mathbb{R}^{d}\right)$ by Theorem 2.8 and $\rho_{s}$ is compactly supported by Theorem 2.5 , we have, for $z \in \partial D_{j}$, that $\rho_{s}(z)=0$. Hence, for any $x \in D_{j}$,
$$
\frac{m \varepsilon}{m-1} \rho_{s}^{m-1}(x)+\left(W * \rho_{s}\right)(x)=\left(W * \rho_{s}\right)(z)
$$
Furthermore, by $(W 4)$, we have that $\left(W * \rho_{s}\right)(z)=\int_{\mathbb{R}^{d}} W(z-y) \rho_{s}(y) d y \leq 0$. Therefore, for all $x \in D_{j}$,
$$
\frac{m \varepsilon}{m-1} \rho_{s}^{m-1}(x) \leq-\left(W * \rho_{s}\right)(x) \leq\|W\|_{L^{1}\left(\mathbb{R}^{d}\right)}\left\|\rho_{s}\right\|_{L^{\infty}\left(\mathbb{R}^{d}\right)},
$$
which implies that
\begin{equation*}
\rho_{s}^{m-1}(x) \leq \frac{m-1}{m \varepsilon}\|W\|_{L^{1}\left(\mathbb{R}^{d}\right)}\left\|\rho_{s}\right\|_{L^{\infty}\left(\mathbb{R}^{d}\right)} \tag{25}
\end{equation*}
for all $x \in D_{j}$. Since $D_{j}$ is an arbitrary connected component of $\operatorname{supp} \rho_{s}$ and the upper-bound of $\rho_{s}^{m-1}$ is independent of $D_{j}$, then Inequality (25) holds for any $D_{j} \subseteq \operatorname{supp} \rho_{s}$ and, thus, for any $x \in \mathbb{R}^{d}$. This yields
$$
\left\|\rho_{s}\right\|_{L^{\infty}\left(\mathbb{R}^{d}\right)} \leq\left(\frac{m-1}{m \varepsilon}\|W\|_{L^{1}\left(\mathbb{R}^{d}\right)}\right)^{1 /(m-2)}
$$
as required. 
It is quite interesting that this boundedness property can be extended to the case of $\omega^{\prime}$ strictly positive on $(0, \infty)$. More precisely, the following theorem holds true.

Theorem 3.2. Let $m>2$. If $\rho_{s}$ is a stationary state of (1) with $W$ satisfying assumptions (W1)-(W3), then
$$
\rho_{s}(x) \leq \rho_{s}^{*}:=\left(\frac{m-1}{m \varepsilon}\|W\|_{L^{1}\left(\mathbb{R}^{d}\right)}\right)^{1 /(m-2)} \text { for all } x \in \mathbb{R}^{d}
$$

Proof. Since $\rho_{s}$ is a stationary state of (1) and is radially decreasing up to a translation on $\mathbb{R}^{d}$, by $[17$, Theorem 2.2$]$, we have that
$$
\frac{m \varepsilon}{m-1} \rho_{s}^{m-1}(x)+\left(W * \rho_{s}\right)(x)=C=\mathrm{const}, \text { for all } x \in \operatorname{supp} \rho_{s}
$$
Since $\rho_{s} \in C\left(\mathbb{R}^{d}\right)$ by [17, Lemma 2.3] and $\rho_{s}$ is compactly supported by [20, Lemma 3.2], we have, for $z \in \partial\left(\operatorname{supp} \rho_{s}\right)$, that $\rho_{s}(z)=0$. Hence,
$$
\frac{m \varepsilon}{m-1} \rho_{s}^{m-1}(x)+\left(W * \rho_{s}\right)(x)=\left(W * \rho_{s}\right)(z)<0
$$
This follows from the fact that $\omega$ is strictly increasing on $(0, \infty)$ and $\lim _{r \rightarrow \infty} \omega(r)=$ 0 . Therefore, $\frac{m \varepsilon}{m-1} \rho_{s}^{m-1}(x)<-\left(W * \rho_{s}\right)(x) \leq\|W\|_{L^{1}\left(\mathbb{R}^{d}\right)}\left\|\rho_{s}\right\|_{L^{\infty}\left(\mathbb{R}^{d}\right)}$, for all $x \in \operatorname{supp} \rho_{s}$ and so
$$
\left\|\rho_{s}\right\|_{L^{\infty}\left(\mathbb{R}^{d}\right)} \leq\left(\frac{m-1}{m \varepsilon}\|W\|_{L^{1}\left(\mathbb{R}^{d}\right)}\right)^{1 /(m-2)}
$$
as required.

\section*{4. Numerical simulations}
In this section we present numerical simulations for Equation (1) to illustrate the results in Theorems 2.10, 3.1, and 3.2. We consider two examples of attractive kernels, one with bounded and the other with unbounded support. We show that in both cases the numerical solutions converge towards stationary states that are radially decreasing up to a translation and compactly supported.

Furthermore, we demonstrate that when $m>2$ the stationary states have an upper-bound independent of the mass. In particular, we are able to confirm numerically that $m=2$ is a threshold such that for $m>2$ there exists a massindependent upper-bound for the density, while for $m \leq 2$, the maximal density increases without bound as more mass is added to the system. In addition, we see that for $m>2$ and sufficiently large mass, the stationary states are approximately constant in the interior of their support. In this case, we show that, for any $W \in L^{1}\left(\mathbb{R}^{d}\right)$, we can approximate the maximum height of the density by minimizing the energy functional. We consider numerical results in one dimension for both Examples 1 and 2. For the initial data, we use a characteristic function, defined on a symmetric real interval. That is, $\rho_{0}=\chi_{[-a, a]}$, for some $a>0$.

In addition, we provide simulations in two-dimensions for the particular case of $W$ compactly supported. We demonstrate that for this choice of interaction kernel, we can obtain stationary solutions which exhibit pattern formation. For these experiments in $2 D$, we vary the choice of initial data in order to illustrate the range of patterns that can be obtained. For our numerical computations, we set $d x=0.4$ and $d t=d x$, where $d x$ and $d t$ denote the sizes of the spatial and time steps, respectively. Unless specified otherwise, for all numerical results, we have set $\varepsilon=1$. Consider the following two choices of attractive kernels:

Example 1. Consider
$$
W(x)=\left\{\begin{array}{l}
-5 e^{1 /\left(|x|^{2}-1\right)},|x|<1 \\
0,|x| \geq 1
\end{array}\right.
$$

Example 2. Consider
$$
W(x)=-e^{-|x|}, x \in \mathbb{R}
$$

We can clearly see from Figures 4, 5, and 6 that, for both Examples 1 and 2, the existence of a mass-independent upper-bound is dependent on the value of the diffusion exponent $m$. More precisely, we see in Figure 4 that, regardless of the size of $\operatorname{supp} W$, for $m=2.1$ and sufficiently large mass, the maximal density of the stationary state remains constant as the mass of the population is increased. That is, the maximal density remains constant with only the support increasing as the mass of the initial data is increased.

In contrast, in Figures 5 and 6, we see that, for $m=2$ and both Examples 1 and 2 , the maximal density continues to grow with the mass of the initial data, implying the absence of a mass-independent upper-bound. Using the same mass as in Figure 4, we obtain convergence towards stationary solutions that do not reach a plateau where the internal density is approximately constant, as is the case for $m>2$. Furthermore, for stationary states with mass that is well above that of the stationary states depicted in Figure 4, the maximum height of the density continues to grow with the mass, suggesting that the maximum height is dependent on the mass for any $M>0$.

Figure 7 shows how stationary states are dependent on the size of $\operatorname{supp} W$. Indeed, we see that, for initial data with sufficiently large support, the support of the stationary solution is made up of multiple connected components when $W$ is given by Example 1 . In contrast, using the same initial data, but where $W$ is given by Example 2, we obtain a stationary state whose support consists of a single component. This is expected for an attractive kernel with infinite support.

\includegraphics[max width=\textwidth, center]{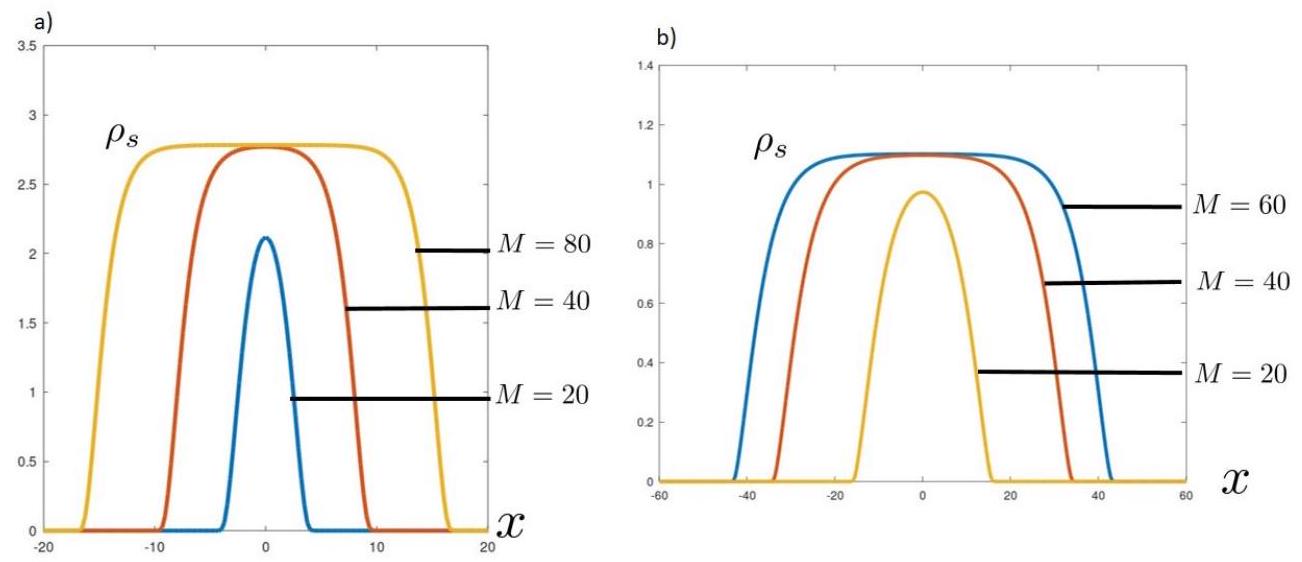}\\
Figure 4: Stationary solutions of Equation (1), where $m=2.1$. a) Example 1: For $M \geq 40$ we have $\max \rho_{s} \approx 2.443$. b) Example 2: For $M \geq 40$ we have $\max \rho_{s} \approx 1.1$

\includegraphics[max width=\textwidth, center]{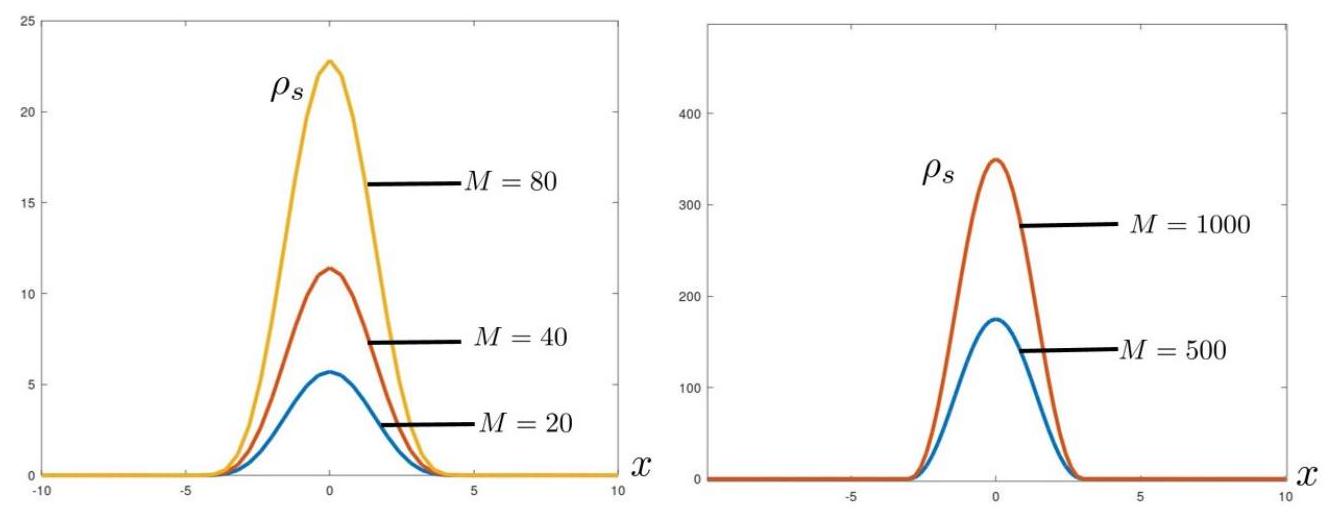}\\
Figure 5: Example 1: Stationary solutions of Equation (1), where $m=2$. The maximum density continues to increase with the mass, while the support remains constant.

\includegraphics[max width=\textwidth, center]{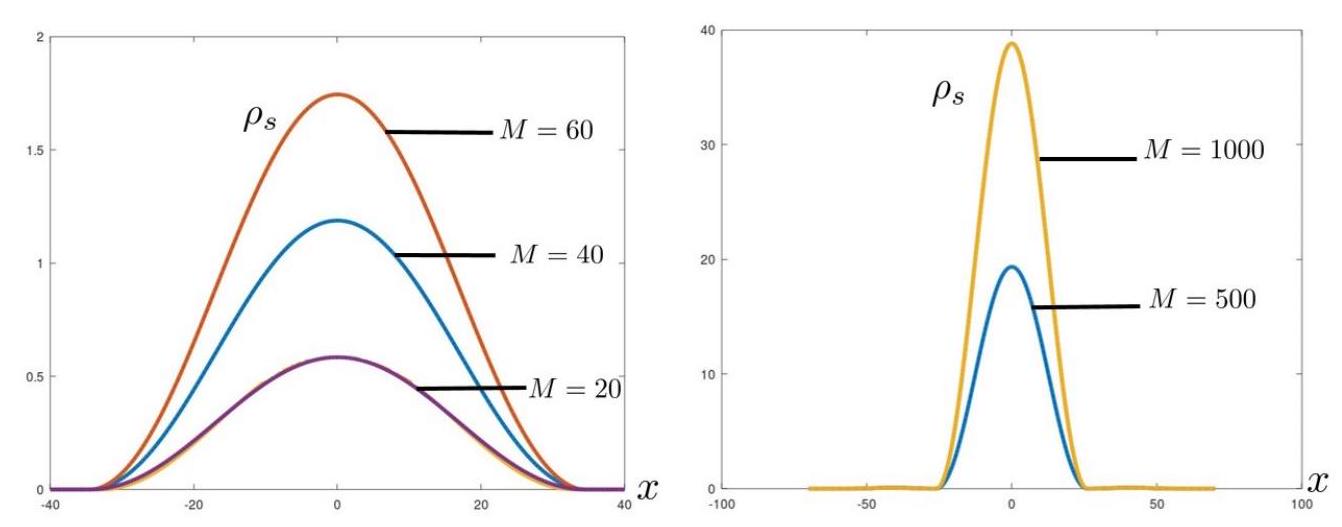}\\
Figure 6: Example 2: Stationary solutions of Equation (1), where $m=2$. The maximum density continues to increase with the mass, while the support remains constant.

\includegraphics[max width=\textwidth, center]{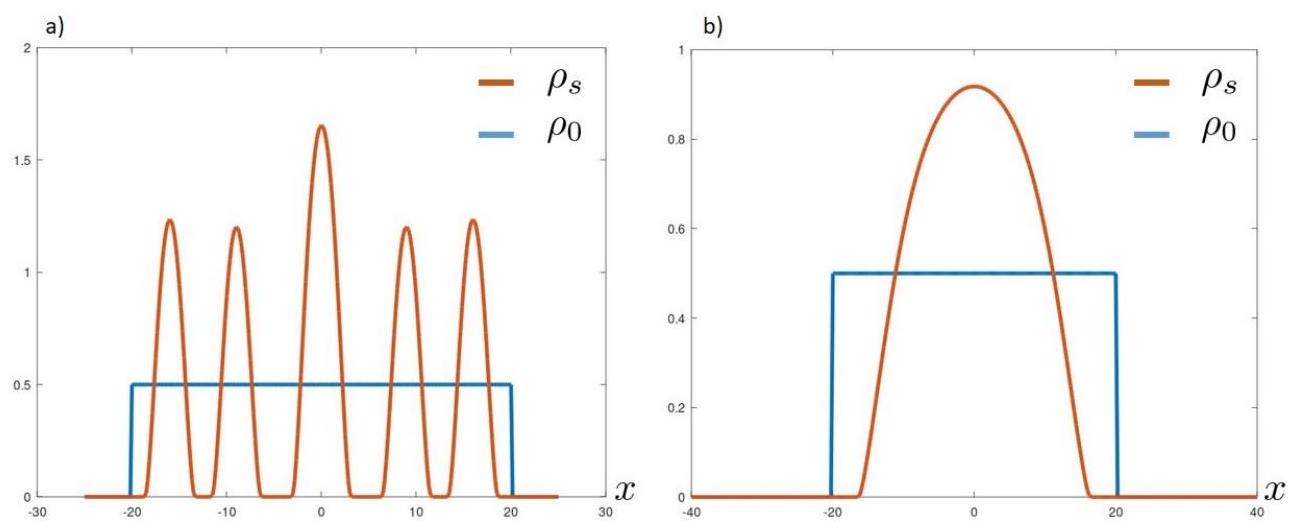}\\
Figure 7: Illustration of the dependence of the stationary solutions on the support of the attractive kernel. a) Stationary solution of Equation (1), where $m=2.1$ and $W$ is given by Example 1, depicting the formation of multiple connected components. b) Stationary solution of Equation (1), where $m=2.1$ and $W$ is given by Example 1. The support of $\rho_{s}$ consists of a single component, contrasting the stationary state with the same initial data given in (a).\vspace{3mm}

\includegraphics[max width=\textwidth, center]{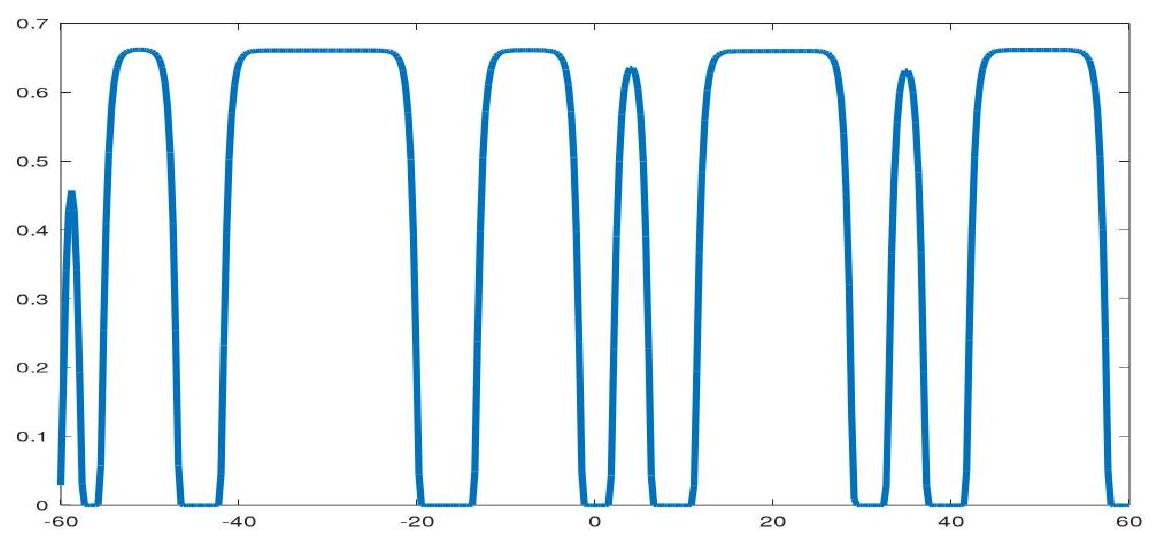}\\
Figure 8: Stationary solution of Equation (1), where $m=3$ and $W(x)=-\max \left(1-|x|^{2}, 0\right)$. The initial condition is randomly distributed with mass $M=50$. The mass of the initial condition is large enough to allow for some swarms to reach their preferred maximum density, where their interior is approximately constant.\vspace{3mm}

\includegraphics[max width=\textwidth, center]{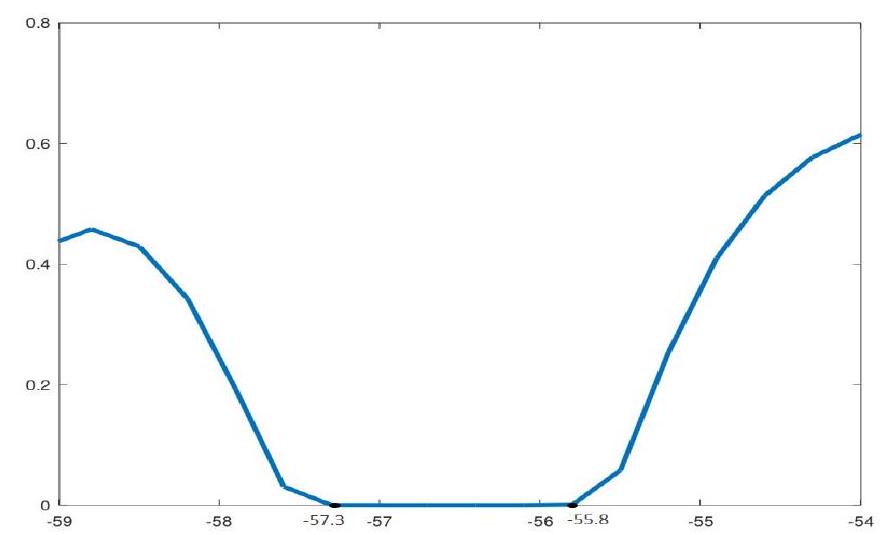}\\
Figure 9: Close up of the two closest swarms of the stationary state depicted in Figure 8, that is, the first and second swarms from the left. The distance between these two swarms is not smaller than the radius of $\operatorname{supp} W$, confirming the analytical result given in Theorem 2.10.2.\vspace{3mm}

Figure 8 gives another example of a stationary state with multiple connected components, for $W$ compactly supported. We choose a randomly distributed initial condition with mass large enough such that some of the swarms in thestationary state reach their preferred maximum density with approximately constant interior. In Figure 9, we provide a close up view of the two swarms in Figure 8 with the shortest proximity, that is, the first and second swarms from the left. We see that the distance between these two swarms is no smaller than the radius of $\operatorname{supp} W$, confirming the analytical result given in Theorem 2.10.2. Furthermore, in all figures for Examples 1 and 2, the stationary states are radially decreasing up to a translation on each connected component of their support, agreeing with the analytical result given in Section 2.

We note that attractive kernels with compact support may be considered advantageous over kernels with infinite support, since they allow for the formation of patterns. This is illustrated in 2D in Figures 10 and 11. We refer also to [9] for video representations of pattern formation for Equation (1) with $W$ compactly supported. In addition, compactly supported interaction kernels can be considered as more realistic with regard to modelling physical and biological aggregations, as the sensing mechanisms of agents in nature do not have an infinite range.

Finally, we see that, in both examples, our numerical results satisfy the properties of conservation of mass and the preservation of positivity of the solution, as expected. Furthermore, when we make the domain large enough to accommodate the mass of the initial condition, we observe that the centre of mass is conserved as the solution evolves. In particular, for the case of $W$ with unbounded support, the centre of mass acts as a centre for the radially decreasing stationary solution.

We note that our numerical simulations are in agreement with the numerical results obtained in [30], where stationary states corresponding to a special case of Equation (1), where $m=3$, are considered. In particular, the authors derive an approximation for the stationary states using the energy functional (3) for the case when the mass is large enough such that the stationary states are approximately constant inside their support. We extend this approximation to stationary states of Equation (1) with diffusion coefficient $m>2$. We consider the energy functional (3) and assume that the mass is large enough such that the stationary states are approximately constant inside their support. Then,

$$
\begin{aligned}
\mathcal{E}\left[\rho_{s}\right] & =\frac{\varepsilon}{m-1} \int_{\mathbb{R}^{d}} \rho_{s}^{m} d x+\frac{1}{2} \int_{\mathbb{R}^{d}} \rho_{s}\left(W * \rho_{s}\right) d x \\
& \approx \frac{\varepsilon}{m-1} \rho_{s}^{m}\left|\operatorname{supp} \rho_{s}\right|_{d}-\frac{\rho_{s}^{2}}{2}|| W \|_{L^{1}\left(\mathbb{R}^{d}\right)}\left|\operatorname{supp} \rho_{s}\right|_{d}
\end{aligned}
$$

Since we assume $\rho_{s}$ is approximately constant in its support, we have that

$$
M=\int_{\operatorname{supp} \rho_{s}} \rho_{s} d x \approx \rho_{s}\left|\operatorname{supp} \rho_{s}\right|_{d}
$$

\includegraphics[max width=\textwidth, center]{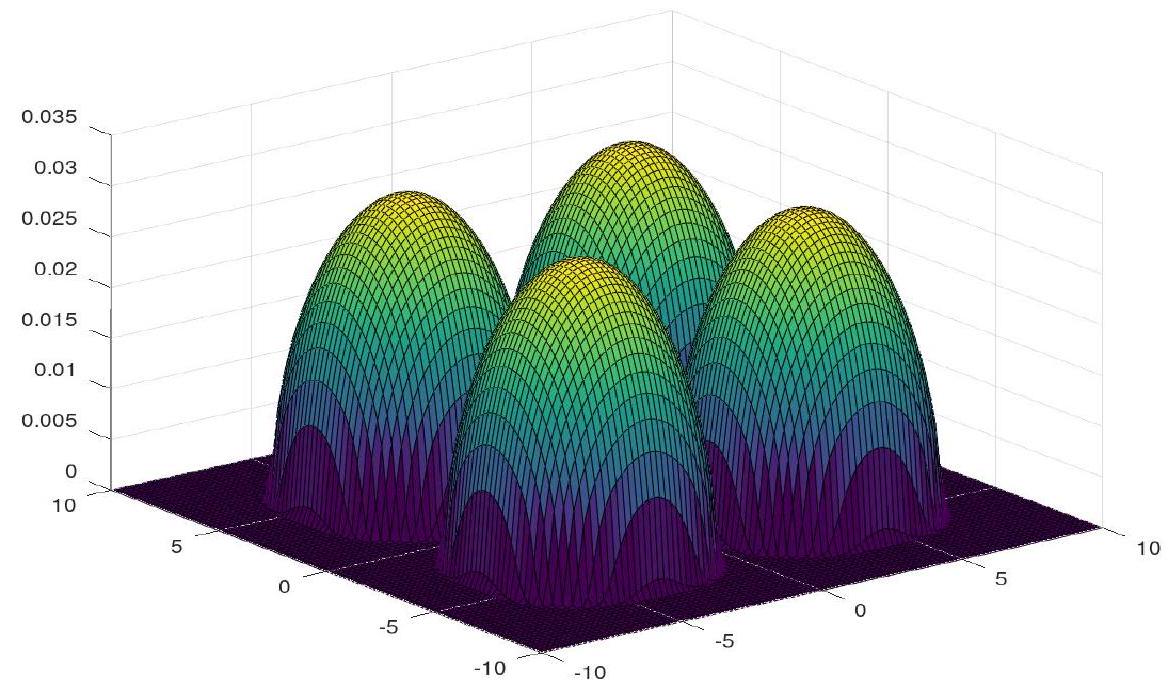}\\
Figure 10: Example 1: Pattern formation of a stationary solution of Equation (1), where $m=3$ and where the initial condition is made up of a linear combination of step functions of the same height that are evenly distributed throughout the domain.

\includegraphics[max width=\textwidth, center]{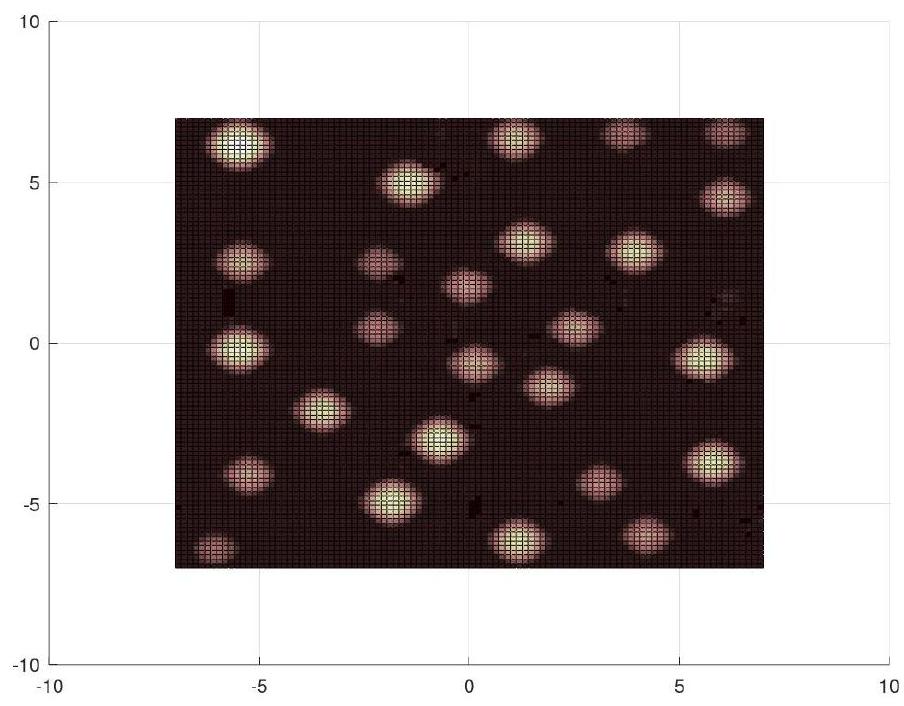}\\
Figure 11: Example 1: Pattern formation of a stationary solution of Equation (1), where $m=3$ and where the initial condition is randomly distributed.

Hence, $\left|\operatorname{supp} \rho_{s}\right|_{d} \approx \frac{M}{\rho_{s}}$, and so

$$
\mathcal{E}\left[\rho_{s}\right] \approx \frac{\varepsilon}{m-1} \rho_{s}^{m-1} M-\frac{\rho_{s}}{2}\|W\|_{L^{1}\left(\mathbb{R}^{d}\right)} M
$$

Now, since $\rho_{s}$ is a stationary state of (1), it is a stationary point of the energy functional, by Theorem 2.9. This yields

$$
\rho_{s} \approx\left(\frac{\|W\|_{L^{1}\left(\mathbb{R}^{d}\right)}}{2 \varepsilon}\right)^{\frac{1}{m-2}}=: \rho_{E}
$$

in $\operatorname{supp} \rho_{s}$.\\
In Table 1, a summary of the results obtained in this section regarding the mass-independent upper-bound of stationary states of (1) is given. We note that the numerical solution for $m=2.1$ and $m=2.5$ goes above $\rho_{E}$, implying that $\rho_{E}$ is not an upper-bound for stationary states of (1). However, we conjecture that this is a numerical error as the numerical stationary state moves above $\rho_{E}$ when approaching the threshold value $m=2$.

Table 1: Approximation of the maximal value of $\rho_{s}$ for Equation (1).

\begin{center}
\begin{tabular}{|c|c|c|c|}
\hline
\multicolumn{4}{|c|}{Example 1} \\
\hline
 & $\rho_{s}^{*}$ & \begin{tabular}{l}
numerical \\
$\max \rho_{s}$ \\
\end{tabular} & $\rho_{E}$ \\
\hline
$\mathrm{m}=2.1$ & 3.942 & 2.443 & 2.476 \\
\hline
$\mathrm{m}=2.5$ & 1.726 & 1.19 & 1.2 \\
\hline
$\mathrm{m}=3$ & 1.333 & 1.0925 & 1.0949 \\
\hline
$\mathrm{m}=3.5$ & 1.333 & 1.0589 & 1.0623 \\
\hline
\multicolumn{4}{|c|}{Example 2} \\
\hline
 & $\rho_{s}^{*}$ & \begin{tabular}{l}
numerical \\
$\max \rho_{s}$ \\
\end{tabular} & $\rho_{E}$ \\
\hline
$\mathrm{m}=2.1$ & 1.59 & 1.1 & 1 \\
\hline
$\mathrm{m}=2.5$ & 1.44 & 1.0117 & 1 \\
\hline
$\mathrm{m}=3$ & 1.333 & 0.99933 & 1 \\
\hline
$\mathrm{m}=3.5$ & 1.27 & 0.99742 & 1 \\
\hline
\end{tabular}
\end{center}

For all numerical results in one dimension, an Euler implicit-explicit (IMEX) scheme, presented in [14], is used. This scheme is based on the explicit scheme devised in [16]. In two dimensions, we again use the Euler IMEX scheme presented in [14], as well as the Alternating-Direction implicit (ADI) method to handle the implicit in time diffusive term [24]. The resulting nonlinear problem is solved using Newton's method. In all simulations, periodic boundary conditions are used.

\section*{5. Conclusion}
In this work, we extend the theory of stationary solutions of Equation (1) to compactly supported attractive kernels. We prove that for $m>2$ the support of any stationary solution is made up of a union of compact, connected components, where the distance between any two connected components is at least the radius of the support of the attractive kernel, if the support is made up of more than one connected component. Furthermore, these stationary solutions are radially decreasing up to a translation on each connected component of their support.

The property of a mass-independent upper-bound of the stationary solutions is relevant in both physical and biological models. Quite remarkably, we show that $m=2$ is a threshold value for both attractive kernels with compact support and infinite support. More precisely, for $m>2$ and for no restriction on the support of the attractive kernel, we prove that stationary solutions have an upper-bound independent of the mass of the solution.

The numerical simulations presented in Section 4 illustrate the analytical results above. We note that our analytical and numerical results agree with characteristics of physical and biological aggregations $[28,29,30]$. That is, the population aggregates to form a group whose spatial extent is bounded. The density also has a mass- independent bound and, for sufficiently large mass, the internal density of the population is approximately constant, implying a preferred inter-organism spacing that is independent of the mass of the population. More precisely, for the case of a domain that is large enough to accommodate the mass, increasing the population beyond a certain point does not cause overcrowding, but rather an increase in the support of the population.

Furthermore, in this work we note that there is an advantage to consider attractive kernels with compact support, as they allow for the formation of patterns, where multiple disjoint "clumps" are formed, in comparison to a single connected "clump" as for the case of an attractive kernel with infinite support.

There are several lines of research arising from this work which can be pursued. For instance, an extension of the current theory can be made regarding Equation (1) defined on a bounded domain. Additionally, since the uniqueness result for stationary states for $\omega^{\prime}$ with unbounded support relies on the radial symmetry property, one may expect that the result on radial symmetry given here may open a window to obtain uniqueness of stationary states on each connected subset of the support for the case of $\omega^{\prime}$ with bounded support. Furthermore, a possibility for future research may be to extend the current theory to the case where both attractive and repulsive terms are modelled non-locally.

\section*{6. Acknowledgements}
This research is partially supported by the DST/NRF SARChI Chair in Mathematical Models and Methods in Bioengineering and Biosciences. The authors would like to further acknowledge the contributions of José A. Carrillo (University of Oxford) and Jacek Banasiak (University of Pretoria) and thank them for their useful comments on this research.

\section*{References}
\begin{flushleft}
[1] M. Adioui O.Arino N.El Saadi, A nonlocal model of phytoplankton aggregation, Nonlinear Analysis: Real World Applications. 6 (2005) 593-607.\\

[2] J. Bedrossian, Intermediate asymptotics for critical and supercritical aggregation equations and Patlak-Keller-Segel models, Commun. Math. Sci. 9 (2011) 1143-1161.\\

[3] J. Bedrossian, Global Minimizers for Free Energies of Subcritical Aggregation Equations with Degenerate Diffusion, Appl. Math. Lett. 24 (2011) 1927-1932.\\

[4] J. Bedrossian, N. Rodriguez, Inhomogeneous Patlak-Keller-Segel models and Aggregation Equations with Nonlinear Diffusion in $\mathbb{R}^{d}$, Discrete Contin. Dyn. Syst. Ser. B. 19 (2014) 1279-1309.\\

[5] J. Bedrossian, N. Rodrguez, A. Bertozzi, Local and Global Well-Posedness for Aggregation Equations and Patlak-Keller-Segel Models with Degenerate Diffusion, Nonlinearity. 24 (2011) 1683-1714.\\

[6] A.J. Bernoff, C.M. Topaz, Biological Aggregation Driven by Social and Environmental Factors: A Nonlocal Model and Its Degenerate Cahn-Hilliard Approximation, SIAM J. Appl. Dyn. Syst. 15 (2016) 15281562.\\

[7] P. Biler, J. Zienkiewicz, Blowing up radial solutions in the minimal KellerSegel model of chemotaxis, J. Evol. Equ. 19 (2019) 71-90.\\

[8] A. Blanchet, J. A. Carrillo, P. Laurenot, Critical mass for a Patlak-KellerSegel model with degenerate diffusion in higher dimensions, Calc. Var. Partial Differ. Equ. 35 (2008) 133-168.\\

[9] C.A. Bright, Nonlinear dynamics (webpage), \href{http://linus.up.ac.za/academic/maths/ND/DAM/index.html}{http://linus.up.ac.za/academic/maths/ND/DAM/index.html}\\

[10] M. Burger, L. Caffarelli, P.A. Markowich, The Keller-Segel model for Chemotaxis with Prevention of Overcrowding: Linear vs. Nonlinear Diffusion, SIAM J. on Math. Analysis . 38 (2014) 1288-1315.\\

[11] M. Burger, V. Capasso, D. Morale, On an aggregation model with long and short range interactions, Nonlinear Analysis: Real World Applications. 8 (2007) 939-958.\\

[12] M. Burger, M. Di Francesco, Y. Dolak-Struss. Partial differential equation models in the socio-economic sciences, Philos Trans A Math Phys Eng Sci. $372(2014) 2028$\\

[13] M. Burger, M. Di Francesco, M. Franek, Stationary states of quadratic diffusion equations with long-range attraction, Communications in Mathematical Sciences . 11 (2013) 709-738.\\

[14] R. Bürger, D. Inzunza, P. Mulet, L.M. Villada, Implicit-explicit methods for a class of nonlinear nonlocal gradient flow equations modelling collective behaviour, Appl. Num. Math. 144 (2019) 234-252.\\

[15] V. Calvez, J. A. Carrillo, F. Hoffmann, Equilibria of homogeneous functionals in the fair-competition regime, Nonlinear Anal. 159 (2017) 85-128.\\

[16] J.A. Carrillo, A. Chertock, Y. Huang, A finite-volume method for nonlinear nonlocal equations with a gradient flow structure, Commun. Comput. Phys. 17 (2014) 233-258.\\

[17] J.A. Carrillo, S. Hittmeir B. Volzone Y. Yao, Nonlinear aggregationdiffusion equations: radial symmetry and long time asymptotics, Inventiones mathematicae. 218 (2019) 889-977.\\

[18] J.A. Carrillo, D. Castorina, B. Volzone, Ground States for Diffusion Dominated Free Energies with Logarithmic Interaction, SIAM J. on Math. Analysis $(2014)$.\\

[19] K. Craig, I. Topaloglu, Aggregation-Diffusion to Constrained Interaction: Minimizers and Gradient Flows in the Slow Diffusion Limit, arXiv preprint arXiv:1806.07415v2 (2019)\\

[20] Delgadino et. al., Uniqueness and non-uniqueness of steady states of aggregation-diffusion equations., arXiv preprint arXiv:1908.09782 (2019)\\

[21] Kaib, G. (2017) Stationary states of an aggregation equation with degenerate diffusion and bounded attractive potential. SIAM J. Math. Anal., 49(1), pp. 272-296.\\

[22] S. Kesavan. Symmetrization and applications. Series in Analysis, 3. World Scientific Publishing Co. Pte. Ltd., Hackensack, NJ, 2006.\\

[23] I. Kim, Y. Yao, The Patlak-Keller-Segel model and its variations: properties of solutions via maximum principle,SIAM J. Math. Anal. 44 (2012) $568-602$.\\

[24] Z. Li, Z. Qiao, T. Tang, Numerical Solution of Differential Equations., Cambridge University Press, 2018.\\

[25] E. H. Lieb, M. Loss, Analysis. Graduate Studies in Mathematics, 14., American Mathematical Society, Providence, RI, 1997.\\

[26] A. Mogilner, L. Edelstein-Keshet. A non-local model for a swarm, J. Math. Biol. 38 (1999) 534570.\\

[27] A. Mogilner, L. Edelstein-Keshet, L. Bent, A. Spiros, Mutual interactions, potentials, and individual distance in a social aggregation, J. Math. Biol. 47 (2003) 353389 .\\

[28] J. K. Parrish, L. Edelstein-Keshet, Complexity, Pattern, and Evolutionary Trade-Offs in Animal Aggregation, Science. 284 (1999) 99-101.\\

[29] J. K. Parrish, W. M. Hamner, Animal Groups in Three Dimensions., Cambridge University Press, Cambridge, UK, 1997.\\

[30] C. M. Topaz, A.L. Bertozzi, M.A. Lewis. A nonlocal continuum model for biological aggregation, Bull. Math. Biol. 68 (2006) 16011623.

\end{flushleft}
\end{document}